\documentclass[12pt]{amsart}

\usepackage{amsthm, amsmath}
\usepackage{amssymb,amsfonts,amscd,verbatim, longtable}
\usepackage[dvips]{graphicx}
\usepackage{xy}
\usepackage{mathrsfs}
\usepackage{wrapfig}

\xyoption{all}
\SelectTips{cm}{10}

\usepackage[hypertex]{hyperref}
\usepackage{multind}
\makeatletter
\renewcommand{\printindex}[2]{\@restonecoltrue\if@twocolumn\@restonecolfalse\fi
  \columnseprule \z@ \columnsep 35pt
  \newpage \twocolumn[{\Large\bf #2 \vskip4ex}]
  \addcontentsline{toc}{chapter}{#2}
  \@input{#1.ind}}
\makeatother

\newtheorem{tr}{Theorem}[section]
\newtheorem*{tr*}{Theorem}

\newtheorem{pr}[tr]{Proposition}
\newtheorem*{pr*}{Proposition}
\newtheorem*{claim}{Claim}
\newtheorem{cor}[tr]{Corollary}
\newtheorem{df}[tr]{Definition}
\newtheorem*{df*}{Definition}
\newtheorem*{not*}{Notation}

\newtheorem*{OQ*}{Conjecture}
\theoremstyle{definition}
\newtheorem{rem}[tr]{Remark}

\def\differential{d}
\renewcommand\d\differential

\newcommand{\ko}{{\mathcal O}}

\newcommand{\ki}{{\mathcal I}}

\newcommand{\E}{{\mathbb E}}

\DeclareMathOperator\im{Im}

\DeclareMathOperator\Ext{Ext}

\DeclareMathOperator\Pic{Pic}

\DeclareMathOperator\id{id}
\DeclareMathOperator\GL{GL}
\DeclareMathOperator\Bl{Bl}
\DeclareMathOperator\Sing{Sing}

\DeclareMathOperator\transpose{T}

\DeclareMathOperator\Supp{{\rm Supp}}

\def\k{\Bbbk}
\renewcommand\Im\im
\def\bb#1{\mathbb #1}
\def\cal#1{\mathcal #1}
\def\fcal#1{\boldsymbol{\cal{#1}}}

\def\ra{\rightarrow}

\def\xra{\xrightarrow}

\def\minus{\smallsetminus}
\def\mto{\mapsto}

\def\pmat#1{\begin{pmatrix}#1\end{pmatrix}}

\def\smat#1{\left(\begin{smallmatrix}#1\end{smallmatrix}\right)}

\def\point#1{\langle #1 \rangle}

\def\refeq#1{$(\ref{#1})$}

\def\P{\bb P}

\def\O{\cal O}

\def\ten{\otimes}
\def\sqten{\boxtimes}

\let\star *

\def\defeq{:=}

\def\iso{\cong}

\def\tilde{\widetilde}

\newtheorem*{refdF}{Definition~\refNUMBER}

\newtheorem*{reftR}{Theorem~\refNUMBER}

\newtheorem*{refcoR}{Corollary~\refNUMBER}

\newtheorem*{refpR}{Proposition~\refNUMBER}

\newlength{\mycolumn}
\title[Modification of the Simpson moduli space $M_{3m+1}(\P_2)$ by vector bundles]{Modification of the Simpson moduli space $M_{3m+1}(\P_2)$ by vector bundles (I)}
\author{Oleksandr Iena}
\address{SISSA\\
via Bonomea, 265\\
34136 Trieste\\
Italy}
\email{iena@sissa.it}
\author{G\"unther Trautmann}
\address{FB Mathematik\\
 Technische Universit\"at Kaiserslautern\\
  Postfach 3049\\
  D-67653 Kai\-sers\-lautern, Germany}
\email{trm@mathematik.uni-kl.de}

\date{}
\begin{document}
\begin{abstract}
We consider the moduli space of  stable vector bundles on curves embedded in $\P_2$ with Hilbert polynomial $3m+1$  and construct a compactification of this space by vector bundles. The result $\tilde M$ is a blow up of the Simpson moduli space $M_{3m+1}(\P_2)$.
\end{abstract}
\maketitle

\tableofcontents
\section{Introduction}

\subsection{Motivation}
Simpson showed in~\cite{Simpson1} that
for an arbitrary smooth projective variety $X$  and for an arbitrary
numerical  polynomial $P\in\bb Q[m]$ there is a coarse moduli space
$M\defeq M_P(X)$ of semi-stable sheaves on $X$ with Hilbert polynomial $P$,
which turns out to be a projective variety.

In many cases $M$ contains an open dense subset $M_B$ whose points
consists of sheaves which are locally free on their support.
So, one could consider $M$ as a compactification of $M_B$. We call sheaves
in the boundary $M\smallsetminus M_B$ \emph{singular}.
It is an interesting question whether and how one could replace the
boundary of singular sheaves by one which consists entirely of vector
bundles with varying and possibly reducible supports.

In this paper we answer this question for the moduli space $M=M_{3m+1}(\P_2)$
of stable sheaves supported on cubic curves in the projective plane
having Hilbert polynomial  $P(m)=3m+1.$ This is the first non-trivial
case of 1-dimensional sheaves on surfaces. It turns out that the blow up
$\tilde{M}$ of $M$ along the locus of singular sheaves is a compactification
in the above sense. Even so this is only a first example it leads to
several interesting constructions which might be helpful in more
general situations.

\subsection{Summary of the paper}
The moduli space $M=M_{3m+1}(\P_2)$ is completely understood
see~\cite{LePotier},~\cite{Freiermuth},~\cite{FreiermuthPreprint}.
The sheaves in $M$ will be called \emph{$(3m+1)$-sheaves}.
Their supports are defined by the their Fitting ideals.
$M$ is smooth of dimension
$10$ and isomorphic to the universal cubic curve. The subvariety $M'$ of
singular sheaves is a smooth subvariety in $M$ of dimension $8$. To
indicate this, we will denote it by $M_8$.

For a singular $(3m+1)$-sheaf $\cal F$ there is only one point
$p\in C'=\Supp\cal F$, where $\cal F$ is not $\ko_{C'}$-free.
This point is a singular point of the curve $C'$.

For each point $p$ of the projective plane, we introduce a reducible surface
$D(p)$ consisting of two irreducible components $D_0(p)$ and $D_1(p)$,
$D_0(p)$ being the blow up of the projective plane at $p$ and $D_1(p)$
being another projective plane, such that these components intersect along the
line $L(p)$ which is the exceptional divisor of $D_0(p)$. Each surface
$D(p)$ can be defined as the subvariety in $\P_2\times\P_2$ with equations
$u_0x_1$, $u_0x_2$, $u_1x_2-u_2x_1$ where the $x_i$ resp. $u_i$ are
the homogeneous coordinates of the first resp. second $\P_2$, such that
the first projection contracts $D_1(p)$ to $p$ and describes $D_0(p)$
as the blow up. We let $\ko_{D(p)}(a,b)$ denote the invertible sheaf
induced by $\ko_{\P_2}(a)\boxtimes\ko_{\P_2}(b)$.

For every singular $(3m+1)$-sheaf we introduce a family of coherent
$1$-dimensional sheaves on $D(p)$, locally free on their (Fitting-)support.
We call the objects of this family $R$-bundles.

\begin{df}\label{df:R-bundles}
An $R$-bundle associated to a singular $(3m+1)$-sheaf $\cal F$
with singular point $p$ is a coherent $1$-dimensional sheaf $\cal E$ on
$D(p)$ subject to the following conditions.

\begin{itemize}
\item $\cal E$ is  locally free on its  support $C=\Supp\cal E$.
\item There is an exact sequence
\begin{equation}\label{eq:resolution back}
0\ra 2\O_{D_1(p)}(-L)\xra{\varrho} \sigma_p^*(\cal F)\xra{\theta} \cal E\ra 0.
\end{equation}
\item $\cal E|_{D_0(p)}$ has Hilbert polynomial $4m+1$ with respect to the
sheaf $\O_{D_0(p)}(1,1)$.
\end{itemize}
\end{df}
$R$-bundles turn out to be flat limits of non-singular $(3m+1)$-sheaves, hence
they can be seen as reasonable replacements for singular $(3m+1)$-sheaves.
$R$-bundles are supported on reducible curves of the type $C=C_0\cup C_1,$
where $C_i=C\cap D_i(p)$.
The curve $C_0$ coincides in most of the cases with the proper transform of
$C'=\Supp\cal F$ under the blow up $D_0(p)\ra \P_2$, in general it is
a subvariety of the total transform of $C'$. The curve $C_1$ is a conic in the
projective plane $D_1(p)$ bearing the degree of an $R$-bundle.
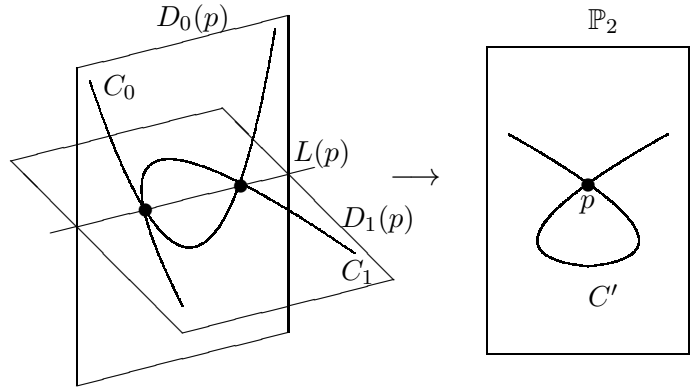
\begin{figure}
\begin{center}
\begin{picture}(250,150)(0,0)
\put(180,20){
\fbox{\begin{picture}(70,110)(-35,-30)
\qbezier(-30,50)(50,5)(0,0)
\qbezier(30,50)(-50,5)(0,0)
\put(0,31){\circle*{5}}
\put(-3,22){\small $p$}
\put(0,-15){\small$C'$}
\put(0,90){\small$\P_2$}
\put(-75,31){$\longrightarrow$}
\end{picture}}
}
\put(0,0){
\begin{picture}
(140,150)(-25,-65)
\put(0,0){\line(4,1){80}}
\put(0,0){\line(0,1){60}}
\put(0,0){\line(0,-1){60}}
\put(80,20){\line(0,-1){60}}
\put(80,20){\line(0,1){60}}
\put(0,60){\line(4,1){80}}
\put(0,-60){\line(4,1){80}}
\put(0,0){\line(1,-1){40}}
\put(80,20){\line(1,-1){40}}
\put(40, -40){\line(4, 1){80}}
\put(0, 0){\line(-1, 1){25}}
\put(80, 20){\line(-1, 1){25}}
\put(-25, 25){\line(4, 1){80}}

\qbezier(5,55)(50,-80)(75,75)
\put(26,6.5){\circle*{5}}
\put(62,15.5){\circle*{5}}
\qbezier(40,-30)(-11,70)(105,-10)
\put(0,0){\line(-4,-1){10}}
\put(80,20){\line(4,1){10}}
\put(82, 25){\small $L(p)$}
\put(10,50){\small ${C_0}$}
\put(100,-20){\small$C_1$}
\put(100,0){\small$D_1(p)$}
\put(30,76){\small$D_0(p)$}
\end{picture}
}
\end{picture}
\end{center}
\caption{Surface $D(p)$ and support of $R$-bundles.}\label{fig:D(p)}
\end{figure}

For $R$-bundles we introduce the following equivalence relation,
which is an ``embedded'' version of Definition~4.1, (ii) from~\cite{Seshadri}
(see also ~\cite{Nagaraj-Seshadri1} and~\cite{Nagaraj-Seshadri2}).
\begin{df}\label{df:equivalence for new sheaves}
Let $\cal E_1$ and $\cal E_2$ be two $R$-bundles on $D(p)$
 associated to the same singular sheaf.
 We call them equivalent if there exists an automorphism $\phi$ of
$D(p)$ that acts identically on $D_0(p)$  and such that
$\phi^*(\cal E_1)\iso \cal E_2$.
\end{df}

The main result of this paper is the following theorem.

\begin{tr}\label{tr:equivalence of sheaves}
The set of equivalence classes of $R$-bundles on $D(p)$ associated to the same
singular $(3m+1)$-sheaf $\cal F$ is in a natural bijection with the
projectivised normal space $\P N_{[\cal F]}$ of the normal space
$N_{[\cal F]} =T_{[\cal F]}(M)/T_{[\cal F]}(M_8)$
to $M_8$ at the point $[\cal F]\in M_8$.
\end{tr}

This justifies the following interpretation.
\begin{cor}
The blow up $\tilde M$ of $M$ along $M_8$  is
the space of all the isomorphism classes of the non-singular $(3m+1)$-sheaves
together with all the equivalence classes of  all $R$-bundles, which are
the points of the exceptional divisor.
\end{cor}

In section \ref{section:universal family} we construct a ``universal''
flat family
of non-singular $(3m+1)$-sheaves and $R$-bundles over $\tilde M$ whose
members are representatives of the equivalence classes of those sheaves.
One can expect that under an appropriate notion of families of
non-singular $(3m+1)$-sheaves together with $R$-bundles the blow up
$\tilde M$ represents the corresponding moduli functor.

\subsection{Structure of the paper}
In Section~\ref{section:space M} we collect the essential facts about the
moduli space $M_{3m+1}(\P_2)$. Details and technicalities necessary to define
$R$-bundles are discussed in  Section~\ref{section:R-bundles}.
In Section~\ref{section:1-fam} we discuss  the most important properties of $R$-bundles.
Using them we prove then the main result in
Section~\ref{section:main result}.
Examples of $R$-bundles are considered in Section~\ref{section:examples}.
 In Section~\ref{section:universal family} we construct  parameter spaces for
all non-singular $(3m+1)$-sheaves (up to isomorphism) and all $R$-bundles
(up to equivalence).

\subsection{Some notations and conventions.}
In this paper $\k$ is an algebraically closed field of characteristic zero.
We work in the category of separated schemes of finite type over $\k$ and
call them varieties, using only their closed points.
Note that we do not restrict ourselves to reduced or irreducible varieties.
Dealing with homomorphism between direct sums of line bundles and identifying
them with matrices, we consider the matrices acting on elements from the
right, i.~e, the composition $X\xra{A}Y\xra{B}Z$ is given by the matrix
$A\cdot B$.
\section{Moduli space $M\defeq M_{3m+1}(\P_2)$}\label{section:space M}
Let us recall here some of the results from~\cite{Freiermuth}.
We consider stable sheaves on $\P_2$ with Hilbert polynomial $3m+1$ and call them simply $(3m+1)$-sheaves. Every $(3m+1)$-sheaf $\cal F$ defines a non-trivial extension
\[
0\ra \O_C\ra \cal F\ra \k_p\ra 0,
\]
where $C$ is a cubic curve supporting $\cal F$ and $\k_p$ is the skyscraper sheaf at $p\in C$ of length $1$, whereby $h^0(\cal F)=1$.

 There exists a fine moduli space $M=M_{3m+1}(\P_2)$ of $(3m+1)$-sheaves.
 $M$ is projective, nonsingular of dimension $10$ and is isomorphic to the universal  cubic
\[
\{(\point{f}, \point{x})\in \P_9\times \P_2\ |\ f(x)=0\},
\] where $\P_9$
 is identified with the space of cubic curves in $\P_2$.
The map underlying this isomorphism is given by $[\cal F]\mto (C, p)$.

The $(3m+1)$-sheaves on $\P_2$ are exactly the sheaves given by locally free resolutions
\[
0\ra 2\O_{\P_2}(-2)\xra{A} \O_{\P_2}(-1)\oplus \O_{\P_2}\ra \cal F\ra 0,
\]
where
\begin{equation}\label{eq:A}
A=\pmat{
z_1&q_1\\
z_2&q_2
}
\end{equation}
with linear independent linear forms $z_1, z_2\in \Gamma(\P_2, \O_{\P_2}(1))$ and non-zero determinant.
The space of all such matrices is a parameter space of $M$ and is denoted by $X$.  $X$ is isomorphic to an open subset in $\k^{18}$ and is acted on by the group
\begin{equation*}
G=\GL_2(\k)\times H,
\end{equation*}
 where $H$ is the group of $2\times 2$ matrices
\[
\pmat{\lambda & z\\0&\mu
},\quad \lambda, \mu\in \k,\quad \lambda\mu\neq 0,\quad z\in \Gamma(\P_2, \O_{\P_2}(1)).
\]
The action is defined  by the rule
\(
(g, h)\cdot A=gAh^{-1}.
\)
$M$ is a geometric quotient of $X$ by $G$, the quotient morphism $X\xra{\nu} M$ is
\[
A=\pmat{z_1&q_1\\z_2&q_2}\mto \point{\det A}\times p(A),
\]
where $p(A)\defeq\point{z_1\wedge z_2}$ denotes the common zero point  of $z_1$ and $z_2$ in $\P_2$.

A $(3m+1)$-sheaf is called singular if it is not locally free on its support. A point $\point{f}\times p\in M$ represents  an isomorphism class of a singular sheaf if and only if $p$ is a singular point of the curve $\{f=0\}\subset \P_2$.
The subvariety of all singular sheaves in $M$ is denoted by $M_8$.
It is easy to verify that $M_8$ is smooth of codimension $2$ in $M$.

The corresponding subvariety in $X$ is denoted by $X_8$.
A matrix $A$ as in~\refeq{eq:A} belongs to $X_8$ if and only if
$q_1(p(A))=q_2(p(A))=0$. These two conditions give two global equations of
$X_8$ in $X$ and one concludes that
$X_8$ is a global complete intersection in $X$, smooth of codimension $2$.

Let $x_0, x_1, x_2\in \Gamma(\P_2, \O_{\P_2}(1))$ be fixed coordinates of $\P_2$. Then  a matrix $A$ from~\refeq{eq:A} with $z_1=x_1$, $z_2=x_2$ belongs to $X_8$ if and only if
\begin{equation}\label{eq:A special}
A=\pmat{
x_1& x_1y_1+x_2y_2\\
x_2&x_1z_1+x_2z_2
}
\end{equation}
for some linear forms $y_i, z_i\in \Gamma(\P_2, \O_{\P_2}(1))$, $i=1,2$.

Let $\tilde M\ra M$ be the blow up of $M$ along $M_8$. Since $M_8$ is smooth of codimension $2$ in $M$, the exceptional divisor $E_M$ of the blow up $\tilde M\ra M$ is isomorphic to the projective normal bundle $\P N_{M_8/M}$.
Let $\tilde X\xra{\alpha} X$ be the blowing up of $X$ along $X_8$.
  Since $X_8$ is defined by two global equations, $\tilde X$ may be considered as a subvariety in $X\times \P_1$ such that the exceptional divisor $E_X$
  of $\tilde X\xra{\alpha} X$ may be identified with $X_8\times \P_1$.

Note that $X_8$ is invariant under the action of $G$. Therefore, since the
blowing up $\alpha:\tilde X\ra X$ is an isomorphism over $X\minus\def\mto{\mapsto} X_8$,
we obtain an action of $G$ on $\tilde X\minus E_X$. This action can be uniquely extended to an action of $G$ on $\tilde X$.  An element $(g, h)\in G$ acts by the rule
\[
(g, h)\cdot (A,\point{t_3, t_4})=(gAh^{-1}, \langle (t_3,
t_4)g^{\transpose} \rangle).
\]
We obtain the following
commutative diagram
\[
\xymatrix{
G\times \tilde X\ar[r]\ar[d]_{\id\times \alpha}&\ar[d]^{\alpha} \tilde
X\\
G\times X\ar[r] & X.
}
\]
Note that for an arbitrary point $(A, \point{t_3, t_4})\in \tilde X$ its stabilizer is the subgroup
\[
St=\left\{
\smat{\lambda&0\\0&\lambda}\times\smat{\lambda&0\\0&\lambda} |\  \lambda\in \k^*
\right\}.
\]
Therefore, we can consider the corresponding free action of the group
$\P G=G/St$ on $\tilde X$.

Note that since $\nu^{-1}(M_8)=X_8$, we obtain a unique lifting $\tilde \nu$ of $\nu$, i.~e., the commutative diagram
\[
\begin{xy}
(0,0)*+{\tilde X}="1";
(20,0)*+{\tilde M}="2";
(0,-15)*+{X}="3";
(20,-15)*+{M.}="4";
{\ar@{->}^{\nu} "3";"4"};
{\ar@{->}^{} "1";"3"};
{\ar@{->}^{} "2";"4"};
{\ar@{->}^{\exists!\ \tilde \nu} "1";"2"};
\end{xy}
\]
Then $\tilde \nu:\tilde X\ra \tilde M$ is $G$-invariant and the set of the
orbits coincides with the set of the fibres $\tilde\nu^{-1}(\xi)$,
$\xi\in \tilde M$. In a  neighbourhood of every point of $\tilde M$ there
is a local section of $\tilde \nu$. Using this and Zariski main theorem
one shows that $\tilde X$ is a principal $\P G$-bundle over $\tilde M$.
Hence $\tilde M$ is a geometrical quotient.

\section{Definition of $R$-bundles and their properties}
\label{section:R-bundles}
\subsection*{Surfaces $D(p)$.}
Let $p$ be a point in $\P_2$.
Let $\sigma_p: \Bl_{0\times p} (\k\times \P_2)\ra \k\times \P_2$ be the blow up of $\k\times \P_2$ at $0\times p$. Consider the composition of $\sigma_p$ with the canonical projection onto $\k$
\[
\Bl_{0\times p} (\k\times \P_2)\xra{\sigma_p} \k\times \P_2\xra{pr_1} \k.
\]
Denote by $D(p)$ the fibre over $0$.
It is a reducible projective surface  consisting of two irreducible components $D_0(p)$ and $D_1(p)$. The first is isomorphic to the blow up $\Bl_p \P_2$ of $\P_2$ at $p$, the second is a projective plane $\P_2$.
Their intersection $L(p)=D_0(p)\cap D_1(p)$ is the exceptional divisor of $D_0(p)$ and is a projective line in $D_1(p)$ (see Figure~\ref{fig:D(p)}).
The restriction $\sigma_p|_{D(p)}$
contracts $D_1(p)$ to $p$ and is the blow up $D_0(p)\ra \P_2$ at $p$.

\subsection*{Invertible sheaves on $D(p)$ and their cohomology.}
Note that all  surfaces $D(p)$ are isomorphic to each other.
Each surface $D(p)$ comes as a closed subvariety of $\P_2\times \P_2$
such that $\sigma_p$ is the restriction of the first projection
to $\P_2$.\hfill\phantom{a}
\begin{wrapfigure}[8]{l}{6.5cm}
\(
\begin{xy}
(0,15)*+{D(p)}="1";
(0,0)*+{\P_2\times \P_2}="2";
(20,0)*+{\P_2}="4";
(-20,0)*+{\P_2}="3";
(-20,15)*+{D_0(p)}="5";
(20,15)*+{D_1(p)}="6";
(-30,0)*+{\{p\}}="7";
{\ar@{>->}"1";"2"};
{\ar@{->}^{\sigma_p}"1";"3"};
{\ar@{->}^-{pr_1}"2";"3"};
{\ar@{->}_-{pr_2}"2";"4"};
{\ar@{->}^{\sigma_p}"5";"3"};
{\ar@{->}_{\iso}"6";"4"};
{\ar@{>->}_{}"6";"1"};
{\ar@{>->}_{}"5";"1"};
{\ar@{->}@`{{**{} ?+/v(-3,2) 130pt /}}_{\sigma_p}"6";"7"};
{\ar@{>->}_{}"7";"3"};
\end{xy}
\)
\end{wrapfigure}
As a subvariety in $\P_2\times \P_2$ every surface $D(p)$ has two different twisting sheaves on $D(p)$, namely
\(
\O_{D(p)}(1, 0)=\O_{\P_2\times \P_2}(1, 0)|_{D(p)}
\), and
\(
\O_{D(p)}(0, 1)=\O_{\P_2\times \P_2}(0, 1)|_{D(p)}.
\)
We can also  define the following two divisors $H$ and $F$ on $D(p)$ by
$\O_{D(p)}(H)=\O_{D(p)}(1,0)$ and $\O_{D(p)}(F)=\O_{D(p)}(0,1)$. In other words
$H$ is defined by the pull-back of a line $h\subset \P_2$ in the
first $\P_2$  and $F$ is
defined by the pull-back of a line $f\subset \P_2$ in the second $\P_2$.

Let $u_0, u_1, u_2$ be the coordinates of the second $\P_2$ and let us choose the coordinates  $x_0, x_1, x_2$  of the first $\P_2$
such that $p=\point{1, 0, 0}$. Then
the surface $D(p)$
is given by the equations
\[
x_1u_2-x_2u_1=0,\quad x_1u_0=0, \quad x_2u_0=0
\]
with $D_0(p)=\{u_0=0\}$, $D_1(p)=\{x_1=x_2=0\}$.
The canonical lifting homomorphisms
\[
\Gamma(\P_2, \O_{\P_2}(h))\ra \Gamma(\P_2, \O_{D(p)}(H)),\quad
\Gamma(\P_2, \O_{\P_2}(f))\ra \Gamma(\P_2, \O_{D(p)}(F))
\]
are isomorphisms. By abuse of notation we denote the images of the homogeneous
coordinates also by
$x_0, x_1, x_2\in \Gamma(\P_2, \O_{D(p)}(H))$ respectively
$u_0, u_1, u_2\in \Gamma(\P_2, \O_{D(p)}(F))$, each forming a basis.
Then we still have the equations $x_1u_2-x_2u_1=0$, $x_1u_0=0$, $x_2u_0=0$ on $D(p)$ with $D_0(p)$, $D_1(p)$ as above.

One can show (cf.~\cite{MyGermanDiss}) that the line bundles $\O_{D(p)}(H)$ and $\O_{D(p)}(F)$ are free generators of the Picard group of $D(p)$. More precisely, the map
\[
\bb Z\oplus \bb Z\ra \Pic(D(p)),\quad (a, b)\mto [\O_{D(P)}(aH+bF)].
\]
is a group isomorphism.

We will denote the restrictions of the divisors $H$ and $F$  to $D_i(p)$, $i=0,1$,
by $H_i$ and $F_i$    respectively. Note that  $H_1\sim 0$, so $H=H_0$.
If it does not cause any misunderstandings, we will often write just $H$ and $F$ for the restrictions $H_i$ and $F_i$.
The intersection line $L(p)$ as a divisor in $D_0(p)$ respectively in $D_1(p)$ will be denoted by $L_0(p)$ respectively $L_1(p)$. There are equivalences of divisors $L_1(p)\sim F_1$ and $L_0(p)\sim H_0-F_0$. The intersection numbers are given by
\begin{align*}
&L_0(p)^2=-1, \quad& &L_1(p)^2=1,\quad& &H_0.F_0=1, \quad& &H_0.L_0(p)=0,\\
&F_0.L_0(p)=1,\quad& &F_1.L_1(p)=1,\quad& &F_0^2=0,\quad& &F_1^2=1.
\end{align*}

\begin{pr}\label{pr:some cohomologies}
(i) The Euler characteristic of the invertible sheaf $\O_{D(p)}(aH+bF)$ is
given by the formula
$\chi(\O_{D(p)}(aH+bF))=\frac{1}{2}(a+b)^2+\frac{3}{2}(a+b)+1.$

(ii) The Hilbert polynomial of $\O_{D(p)}(aH+bF)$ with respect to the
invertible sheaf $\cal L=\O_{D(p)}(H+F)$ equals
\begin{equation}\label{eq:Hilbert D(p)}
2m^2+[2(a+b)+3]\cdot m+\frac{1}{2}(a+b)^2+\frac{3}{2}(a+b)+1.
\end{equation}

(iii) Higher cohomology groups vanish for the following sheaves:
\begin{multline*}
\O_{D(p)}(-2F),\quad
\O_{D(p)}(-H),\quad
\O_{D(p)}(-F),\quad
\O_{D(p)},\quad
\O_{D(p)}(H-F),\quad
\O_{D(p)}(-H+F),\\
\O_{D(p)}(H),\quad
\O_{D(p)}(F),\quad
\O_{D(p)}(-H-F),\quad
\O_{D(p)}(H+F),\quad
\O_{D(p)}(H+2F),
\end{multline*}
hence their $h^0$ can be computed using the formula above.
\end{pr}
\begin{proof}[Some key tools for the proof]
All the statements of Proposition~\ref{pr:some cohomologies} can be directly verified using the following observations.

For a locally free sheaf $\cal G$ on $D(p)$ there is the ``gluing'' exact sequence
\[
0\ra \cal G\ra \cal G|_{D_0(p)}\oplus \cal G|_{D_1(p)}\ra \cal G|_{L(p)}\ra 0
\]
and in particular the exact sequences
\begin{equation}\label{eq:gluing}
0\ra \O_{D(p)}(aH+bF)\ra \O_{D_0(p)}(aH_0+bF_0)\oplus \O_{D_1(p)}(bF_1)\ra \O_L(b)\ra 0, \quad a, b\in \bb Z.
\end{equation}
The cohomology of the sheaves $\O_{D(p)}(aH+bF)$ can be deduced from~\refeq{eq:gluing}
and the exact sequences
\begin{equation*}\label{eq:resolution D_0}
0\ra \O_{\P_2\times \P_1}(a-1,b-1)\ra \O_{\P_2\times \P_1}(a, b) \ra \O_{D_0(p)}(aH_0+bF_0)\ra 0
\end{equation*}
describing the embedding $D_0(p)\subset \P_2\times \P_1$ with equation $u_1x_2-u_2x_1=0$ using the K\"unneth formulas (see~\cite{Sampson}) for $\P_2\times\P_1$.
\end{proof}
 \subsection*{Some canonical homomorphisms.} In the following we need some
canonical homomorphisms related to the reducible surface $D(p)$:
\begin{itemize}
\item There is a canonical section $\O_{D(p)}\xra{s} \O_{D(p)}(H-F)$ induced by the canonical section of the canonical divisor $L_0\sim H-F_0$
via the gluing sequence~\refeq{eq:gluing}, vanishing along $D_1(p)$.

\item For any $a\in \bb Z$ there is the homomorphism
\(
\O_{D_1(p)}\xra{u_0}\ko_{D(p)}(aH+F)
\)
induced by the diagram
\[
\begin{xy}
(0,0)*+{0}="0";
(20,0)*+{\O_{D(p)}(aH+F)}="1";
(70,0)*+{\O_{D_0(p)}(aH_0+F_0)\oplus \O_{D_1(p)}(F_1)}="2";
(110,0)*+{\O_L(1)}="3";
(125,0)*+{0}="4";
(70,15)*+{\O_{D_1(p)}}="5";
{\ar@{->}"0";"1"};
{\ar@{->}"1";"2"};
{\ar@{->}"2";"3"};
{\ar@{->}"3";"4"};
{\ar@{-->}_{u_0}"5";"1"};
{\ar@{->}^{(0, u_0)}"5";"2"};
\end{xy}
\]
because $u_0$ is the equation of $L$ in $D_1(p)$.

\item Using the gluing sequence~\refeq{eq:gluing} one obtains the exact sequence
\[
0\ra \O_{D_1(p)}(-F_1)\xra{u_0}\O_{D(p)}(aH)\xra{s}\O_{D(p)}((a+1)H-F)\xra{r_1} \O_{D_1(p)}(-F_1)\ra 0,
\]
where $r_1$ denotes the restriction homomorphism to $D_1(p)$.

\item The sections $x_1, x_2$ of $\O_{D(p)}(H)$ factorize as $x_\nu=u_\nu\circ s$:
\begin{equation}\label{eq:decomposition of x_i}
\begin{xy}
(0,0)*+{\O_{D(p)}(-H)}="1";
(30,0)*+{\O_{D(p)}.}="2";
(15,-15)*+{\O_{D(p)}(-F)}="3";
{\ar@{->}^-{x_\nu}"1";"2"};
{\ar@{->}_{s}"1";"3"};
{\ar@{->}_{u_\nu}"3";"2"};
\end{xy}
\end{equation}
\end{itemize}

\subsection*{$R$-bundles and their properties.} We define $R$-bundles as in Definition~\ref{df:R-bundles}.
In addition to the three items in the definition, $R$-bundles
have four other properties, which will be derived at the end of Section \ref{section:1-fam}.

\begin{pr}\label{pr:R properties}
Let $\cal E$ be an $R$-bundle and let $C_i=D_i(p)\cap C$, $i=0, 1$, be the components of its support. Then
\begin{enumerate}
\item $\cal E$ has a locally free resolution
\begin{equation}\label{eq:resolution new}
0\ra 2\O_{D(p)}(-H-F)\xra{}\O_{D(p)}(-H)\oplus \O_{D(p)}\ra
\cal E\ra 0.
\end{equation}
\item
The restriction of ~\refeq{eq:resolution new}  to $D_1=D_1(p)$ induces a
resolution
\[
0\ra 2\O_{D_1(p)}(-L) \xra{} 2\O_{D_1(p)}\xra{} \cal E_{D_1(p)}\ra 0
\]
such that $\cal E|_{D_1(p)}$ is a semistable $(2m+2)$-sheaf on
$D_1(p)\iso \P_2$ with $\Supp (\cal E|_{D_1(p)})$ a conic.

\item $\cal E|_{D_0(p)}$ is isomorphic to the structure sheaf $\O_{C_0}$
with the resolution
\[
0\ra \O_{D_0(p)}(-2F_0-H_0)\xra{}\O_{D_0(p)}\xra{} \O_{C_0}\ra 0.
\]

\item $h^0 \cal E=1$, and the non-zero section gives rise to a non-trivial
extension sequence
\[
0\ra \O_{C}\ra \cal E\ra \k_q\ra 0,
\]
where $q\in C_1\minus L$ is uniquely determined by $\cal E$.
\end{enumerate}
\end{pr}

\begin{rem}
Let $\cal E$ be an $R$-bundle. Let  $C'$ be the cubic curve $C'=\Supp \cal F$.
Note that the curve $C_0$ is a subvariety of the total transform of $C'$ under
the blow up $D_0(p)\ra \P_2$. In most of the cases it is the proper transform of $C'$.
Therefore, the restriction map $\sigma_p:C_0\ra C'$
may be considered as a ``partial normalization'' of $C'$.
\end{rem}
\begin{rem}
Note that the properties 1) and 4) of Proposition~\ref{pr:R properties} are analogous to those of the $(3m+1)$-sheaves. The degree $1$ sheaf $\cal F$ on $C'$ is replaced by $\cal E$ whose degree has been shifted to the additional curve $C_1$ with $\cal E_{D_1}$ a $(2m+2)$-sheaf and a vector bundle of degree $1$ on $C_1$.
\end{rem}

\begin{rem}\label{rem:Beilinson D_i}
1) Note that the restriction of resolution~\refeq{eq:resolution new} to
the component $D_1(p) = \P_2$ is a Beilinson resolution of $\cal E|_{D_1(p)}$
on $\P_2$.

2) One can also show that the restriction of~\refeq{eq:resolution new} to $D_0(p)$ is the resolution of Beilinson type, see \cite{Ancona-Ottaviani}, Theorem~8.
\end{rem}
\begin{rem}\label{rem:push forward}
Note that ${\sigma_p}_*\sigma_p^*(\cal F) \iso \cal F$
 for every $(3m+1)$-sheaf $\cal F$. Applying ${\sigma_p}_*$ to
 sequence~\refeq{eq:resolution back} and using that
 \(
 R^0{\sigma_p}_*\cal \O_{D_1(p)}(-L)=R^1{\sigma_p}_*\cal \O_{D_1(p)}(-L)=0
 \)
we obtain the isomorphism ${\sigma_p}_*\cal E\iso \cal F$.
\end{rem}

The proof of Proposition~\ref{pr:R properties} will be a consequence
of the description of $R$-bundles as flat limits of non-singular
$(3m+1)$-sheaves in the following section.

\section{$R$-bundles as $1$-dimensional degenerations of $(3m+1)$-sheaves.}\label{section:1-fam}
Let $A$ be a matrix in $X_8$ and let $B$ be a matrix representing a morphism
$2\O_{\P_2}(-2)\ra \O_{\P_2}(-1)\oplus \O_{\P_2}$. Recall (cf. Section~\ref{section:space M}) that $A$ and $B$ can be considered as elements in $\k^{18}$.
Consider the
morphism
\[
l_B:\k \ra \k^{18},\quad t\mto A+tB.
\]
Let $T=l_B^{-1}(X)$. This way we obtain the morphism
\begin{equation}\label{eq:l_B}
l_B:= l_B|_T:T\ra X.
\end{equation}

By the property of the space $M$ we obtain a $(3m+1)$-family $\fcal F$
over $T$ with the resolution
\begin{equation}\label{eq:basic family}
0\ra 2\O_{T\times\P_2}(-2H)\xra{A+tB} \O_{T\times\P_2}(-H)\oplus
\O_{T\times\P_2}\ra \fcal F\ra 0.
\end{equation}
Here $H$ is represented by the pull-back of a line $h\subset \P_2$.
We choose $h$ such that the point $p$ does not lie on $h$.
By shrinking $T$ we can also assume that
$A+tB\in X\minus X_8$ for all $t\in T$, $t\neq 0$.
In other words, the restrictions
$\fcal F_t$ of the sheaf $\fcal F$ to the fibres $t\times \P_2\iso \P_2$ are
non-singular $(3m+1)$-sheaves  for all $t\in T$, $t\neq 0$. So
the singular $(3m+1)$-sheaf $\fcal F_0$ is a  flat $1$-parameter degeneration
of non-singular $(3m+1)$-sheaves.

Let $p=p(A)$ and
consider the blow up
$Z:= \Bl_{0\times p}(T\times \P_2)\xrightarrow{\sigma}T\times \P_2$.
As above we denote  the exceptional divisor of $\sigma$ by  $D_1=D_1(p)$.
By abuse of notation the lifting of the divisor $H\subset T\times \P_2$ is
again denoted by $H$.
\begin{rem}\label{rem:embz} Letting $x_i$ denote the homogeneous coordinates
of $\P_2$, such that the point $p$ has the equations $tx_0, x_1, x_2,$
$Z$ is embedded in $T\times\P_2\times\P_2$ with equations
\begin{equation*}
tx_0u_1-x_1u_0,\quad tx_0u_2-x_2u_0,\quad x_1u_2-x_2u_1,
\end{equation*}
where the $u_i$ are the coordinates of the second  $\P_2$. It follows
that the fibre $Z_0$ for $t=0$ equals $D(p)$, see section
\ref{section:R-bundles}.
\end{rem}
Note that the morphism $Z\xra{\sigma} T\times \P_2\xra{pr_1} T$ is flat.
Indeed,
since both $Z$ and $T$ are regular, $\dim Z=3$, $\dim T=1$, and $\dim
Z_t=2=\dim Z-\dim T$ for all $t\in T$, this follows
from~\cite{EGA-IV-II}, 6.1.5.
Applying $\sigma^*$ to sequence \refeq{eq:basic family} we obtain
the sequence
\[
0\ra 2\O_{Z}(-2H)\xra{\sigma^*(A+tB)} \O_{Z}(-H)\oplus \O_{Z}\ra
\sigma^*(\fcal F)\ra 0,
\]
which remains exact because the sheaf $O_{Z}(-2H)$ is locally free and,
therefore, has no torsion.
There is a canonical section $s\in\Gamma(Z, \O_Z(D_1))$,
which gives us  the exact sequence
\[
0\ra \O_Z(-D_1)\xra{s} \O_Z\ra \O_{D_1}\ra0.
\]
Tensoring with $\O_Z(D_1-2H)$ one gets the exact
sequence
\[
0\ra 2\O_{Z}(-2H)\xra{\smat{s&0\\0&s}}2\O_{Z}(-2H+D_1)\ra 2\O_{D_1}(-L)\ra 0.
\]
We use here that $H$ and $D_1$ do not meet (our choice of $H$) and that $\O_{D_1}\ten \O_Z(D_1)\iso \O_{D_1}(-L)$ (properties of blow ups).

Note that $A+tB$ vanishes at $0\times p$. Therefore, the morphism  $\sigma^*(A+tB)$ vanishes on $D_1$ and hence
factorizes uniquely through $s$, i. e., there exists
\[
2\O_{Z}(-2H+D_1)\xra{\phi(A, B)}\O_{Z}(-H)\oplus \O_{Z}
\]
 such that the diagram
\[
\begin{xy}
(0,0)*+{2\O_{Z}(-2H)}="1";
(50,0)*+{\O_{Z}(-H)\oplus
\O_{Z}}="2";
(25,-15)*+{2\O_{Z}(-2H+D_1)}="3";
{\ar@{->}^-{\sigma^*(A+tB)}"1";"2"};
{\ar@{->}_-(.3){\smat{s&0\\0&s}}"1";"3"};
{\ar@{->}_-(.6){\phi(A, B)}"3";"2"};
\end{xy}
\] commutes.

Note that $\phi(A, B)$ is injective since $2\O_{Z}(-2H+D_1)$ is torsion
free and since $\smat{s&0\\0&s}$ is an isomorphism outside of the exceptional divisor $D_1$.

Note also that  the exceptional divisor $D_1$ is equivalent to  the difference
$H-F$, where $F$ is the pull-back of a line in the second $\P_2$ along the
standard embedding $Z\subset T\times \P_2\times \P_2$.
Hence we obtain the following commutative diagram with exact rows and columns.
\begin{equation}\label{eq:main construction diagram}
\begin{xy}
(105,20)*+{0}="35";
(40,10)*+{0}="23";
(105,10)*+{2\O_{D_1}(-L)}="25";
(15,0)*+{0}="12";
(40,0)*+{ 2\O_{Z}(-2H)}="13";
(80,0)*+{ \O_{Z}(-H)\oplus \O_{Z}}="14";
(105,0)*+{\sigma^*\fcal F}="15";
(120,0)*+{ 0}="16";
(15,-15)*+{0}="02";
(40,-15)*+{2\O_{Z}(-H-F)}="03";
(80,-15)*+{ \O_{Z}(-H)\oplus \O_{Z}}="04";
(105,-15)*+{\tilde{\fcal F}}="05";
(120,-15)*+{ 0.}="06";
(40,-25)*+{2\O_{D_1}(-L)}="-13";
(105,-25)*+{0}="-15";
(40,-35)*+{0}="-23";
{\ar@{->}"12";"13"};
{\ar@{->}^-{\sigma^*(A+tB)}"13";"14"};
{\ar@{->}"14";"15"};
{\ar@{->}"15";"16"};
{\ar@{->}"02";"03"};
{\ar@{->}^-{\phi(A, B)}"03";"04"};
{\ar@{->}"04";"05"};
{\ar@{->}"05";"06"};
{\ar@{->}"23";"13"};
{\ar@{->}^{\smat{s&0\\0&s}}"13";"03"};
{\ar@{->}"03";"-13"};
{\ar@{->}"-13";"-23"};
{\ar@{->}"35";"25"};
{\ar@{->}"25";"15"};
{\ar@{->}"15";"05"};
{\ar@{->}"05";"-15"};
{\ar@{=}"14";"04"};
\end{xy}
\end{equation}
where $\tilde{\fcal F}$ is defined as cokernel.
\begin{pr}\label{pr:local freeness}
The sheaf $\tilde{\fcal F}$ is locally free on
its support if and only if $B$ is a normal vector to $X_8$ at $A$, i. e., if and only if  $B\in T_A X\minus T_A X_8$.
\end{pr}
\begin{proof}[Sketch of the proof.]
The sheaf $\tilde{\fcal F}$ is not locally free  at some point if and only
if  the morphism  $\phi(A, B)$ vanishes at this point.
Since the only zero point of the morphism $A+tB$ is $(0, p)$ and since the preimage of
$(0, p)$ is the exceptional divisor $D_1$, we  conclude that $\phi(A, B)$
may only vanish at points lying in  $D_1$.

To simplify the considerations one can assume without loss of generality that $p=\point{1,0,0}$ and that $A$ is of the form~\refeq{eq:A special}.
Let
\begin{equation}\label{eq:A special coeff}
A=\pmat{
x_1& A_{01}x_0x_1+\dots+A_{22}x_2^2\\
x_2&B_{01}x_0x_1+\dots+B_{22}x_2^2
},\quad
B=
\pmat{
\xi_0x_0+\xi_1x_1+\xi_2x_2&\xi_{00}x_0^2+\dots+\xi_{22}x_{2}^2\\
\eta_0x_0+\eta_1x_1+\eta_2x_2&\eta_{00}x_0^2+\dots+\eta_{22}x_{2}^2
},
\end{equation} then
straightforward calculations show that vanishing of $\phi(A, B)$ at a point $q\in D_1(p)$ is equivalent to the system
\begin{equation}
\label{eq:tangent equations}
\begin{cases}
\xi_{00}=A_{01}\xi_0+A_{02}\eta_0,\\
\eta_{00}=B_{01}\xi_0+B_{02}\eta_0.
\end{cases}
\end{equation}
One easily checks that~\refeq{eq:tangent equations} are
tangent equations at $A$.
This  proves the required statement.
\end{proof}

The sheaf $\tilde {\fcal F}$ is flat over $T$, because its resolution
remains exact after restriction to fibres.
So we obtain for every normal vector $B$ a sheaf
\begin{equation}\label{eq:E(A, B)}
\cal E=\cal E(A, B)\defeq \tilde{\fcal F}|_{D(p)}
\end{equation}
on $D(p)$.
This sheaf is locally free on its support and is a flat 1-dimensional
degeneration of non-singular $(3m+1)$-sheaves.
Using flatness of $\tilde{\fcal F}$ and restricting
diagram~\refeq{eq:main construction diagram} to $D(p)$
one obtains an exact sequence
\[
0\ra 2\O_{D_1(p)}(-L)\ra \sigma_p^*(\fcal F_0)\ra \cal E\ra 0
\]
and the locally free resolution of $\cal E$ on $D(p)$,
\begin{equation}\label{eq:Phi(A, B)}
0\ra 2\O_{D(p)}(-H-F)\xra{\Phi(A, B)}  \O_{D(p)}(-H)\oplus \O_{D(p)}\ra \cal E\ra 0, \quad \Phi(A, B)\defeq \phi(A, B)|_{D(p)}.
\end{equation}
Straightforward calculations show now that the sheaf $\cal E$ satisfies
Definition~\ref{df:R-bundles}.

We have obtained a construction that produces $R$-bundles.
The following proposition shows that every $R$-bundle can be obtained by this
construction for some $A$ and $B$ and yields at the same time a proof of
Proposition~\ref{pr:R properties}.

\begin{pr}\label{pr:all R-bundles are degenerations}
Each $R$-bundle $\cal E$ is part of an exact diagram
\begin{equation}\label{eq:diagram}
\begin{xy}
(0,25)*+{}="31";
(15,25)*+{}="32";
(40,25)*+{ 0}="33";
(110,25)*+{0}="35";
(125,25)*+{}="36";
(0,15)*+{}="21";
(15,15)*+{}="22";
(40,17)*+{ 2\O_{D_1}(-L)}="23";
(80,15)*+{}="24";
(110,15)*+{2\O_{D_1}(-L)}="25";
(125,15)*+{}="26";
(-8,0)*+{0}="11";
(8,0)*+{2\O_{D_1}(-L)}="12";
(40,0)*+{ 2\O_{D(p)}(-2H)}="13";
(80,0)*+{ \O_{D(p)}(-H)\oplus\O_{D(p)}}="14";
(110,0)*+{\sigma^*\cal F }="15";
(125,0)*+{ 0}="16";
(0,-15)*+{}="01";
(15,-15)*+{0}="02";
(40,-15)*+{2\O_{D(p)}(-H-F)}="03";
(80,-15)*+{ \O_{D(p)}(-H)\oplus\O_{D(p)}}="04";
(110,-15)*+{\cal E}="05";
(125,-15)*+{ 0}="06";
(0,-25)*+{}="-11";
(15,-25)*+{}="-12";
(40,-25)*+{2\O_{D_1}(-L)}="-13";
(75,-25)*+{}="-14";
(110,-25)*+{0}="-15";
(125,-25)*+{}="-16";
(40,-35)*+{0}="-23";
{\ar@{->}"11";"12"};
{\ar@{->}_-{\smat{u_0&0\\0&u_0}}"12";"13"};
{\ar@{->}^-{\sigma^*A}"13";"14"};
{\ar@{->}^-{\pi}"14";"15"};
{\ar@{->}"15";"16"};
{\ar@{->}"02";"03"};
{\ar@{->}^-{\Phi}"03";"04"};
{\ar@{->}^-{\pi'}"04";"05"};
{\ar@{->}"05";"06"};
{\ar@{->}"33";"23"};
{\ar@{->}^{\smat{u_0&0\\0&u_0}}"23";"13"};
{\ar@{->}^{\smat{s&0\\0&s}}"13";"03"};
{\ar@{->}"03";"-13"};
{\ar@{->}"-13";"-23"};
{\ar@{->}"35";"25"};
{\ar@{->}^{\varrho}"25";"15"};
{\ar@{->}^{\theta}"15";"05"};
{\ar@{->}"05";"-15"};
{\ar@{=}"14";"04"};
{\ar@{=}"12";"23"};
\end{xy}
\end{equation}
with $\Phi=\Phi(A, B)$ for some $A\in X_8$ and
$B\in T_A X\minus T_A X_8$.

If $p=\point{1, 0, 0}$ and $A$ is as in~\refeq{eq:A special}, then
\[
\Phi=\pmat{
u_1&u_1y_1+u_2y_2\\
u_2&u_1z_1+u_2z_2
}+
\pmat{\xi_0&\xi_{00}x_0\\\eta_0&\eta_{00}x_0}u_0.
\]
\end{pr}
\begin{proof}
We divide the proof into the following steps.

1) Let $\cal E$ be an $R$-bundle as in Definition~\ref{df:R-bundles}.
 For the proof we may assume that $\cal F$ is the cokernel of an $A$ as in~\refeq{eq:A special}. Let $\sigma=\sigma_p$, $D=D(p)$, $D_i=D_i(p)$. Then $\sigma^*(\cal F)$ has the resolution
\[
0\ra 2\O_{D_1(-L)}\xra{\smat{u_0&0\\0&u_0}}  2\O_{D}(-2H)\xra{\sigma^*A} \O_{D}(-H)\oplus\O_{D}\xra{\pi}
\sigma^*\cal F\ra 0.
\]
\par
\noindent
2) The homomorphism $\varrho$ can be uniquely lifted to a morphism of resolutions
\begin{equation}\label{eq:lifting res}
\begin{split}
\begin{xy}
(0,0)*+{2\O_{D}(-H-F)}="1";
(30,0)*+{2\O_{D_1}(-L)}="2";
(0,-15)*+{\O_{D}(-H)\oplus \O_{D}}="3";
(30,-15)*+{\sigma^*\cal F}="4";
{\ar@{->}"1";"2"};
{\ar@{->}^-{\pi}"3";"4"};
{\ar@{->}_{\tilde A}"1";"3"};
{\ar@{->}^{\varrho}"2";"4"};
(-35,0)*+{2\O_{D}(-2H)}="5";
(-35,-15)*+{2\O_{D}(-2H)}="6";
{\ar@{->}^{B}"5";"6"};
{\ar@{->}^-{\sigma^*A}"6";"3"};
{\ar@{->}^-{\smat{s&0\\0&s}}"5";"1"};
(-65,0)*+{2\O_{D_1}(-L)}="7";
(-65,-15)*+{2\O_{D_1}(-L)}="8";
{\ar@{->}^-{\smat{u_0&0\\0&u_0}}"8";"6"};
(45,0)*+{0}="9";
(45,-15)*+{0}="10";
{\ar@{->}"2";"9"};
{\ar@{->}"4";"10"};
{\ar@{>}^{B}"7";"8"};
{\ar@{->}^-{\smat{u_0&0\\0&u_0}}"7";"5"};
(-80,0)*+{0}="11";
(-80,-15)*+{0}="12";
{\ar@{->}"11";"7"};
{\ar@{->}"12";"8"};
\end{xy}
\end{split}
\end{equation}
because of the vanishing of the relevant $\Ext$-groups, following from Proposition~\ref{pr:some cohomologies}.
\begin{claim}
$B$ is an isomorphism.
\end{claim}
\begin{proof}[Proof of the Claim.]
The restriction of~\refeq{eq:lifting res} to $D_0$ becomes the exact diagram
\[
\begin{xy}
(0,0)*+{2\O_{D_0}(-H-F_0)}="1";
(30,0)*+{2\O_{L}(-1)}="2";
(0,-15)*+{\O_{D_0}(-H)\oplus \O_{D_0}}="3";
(30,-15)*+{\sigma_0^*\cal F}="4";
{\ar@{->}"1";"2"};
{\ar@{->}^-{\pi}"3";"4"};
{\ar@{->}_{\tilde A_{D_0}}"1";"3"};
{\ar@{->}^{\varrho_{D_0}}"2";"4"};
(-35,0)*+{2\O_{D_0}(-2H)}="5";
(-35,-15)*+{2\O_{D_0}(-2H)}="6";
{\ar@{->}^{B}"5";"6"};
{\ar@{->}^-{\sigma_0^*A}"6";"3"};
{\ar@{->}^-{\smat{s_0&0\\0&s_0}}"5";"1"};
(-55,0)*+{0}="7";
(-55,-15)*+{0}="8";
{\ar@{->}"8";"6"};
(45,0)*+{0}="9";
(45,-15)*+{0,}="10";
{\ar@{->}"2";"9"};
{\ar@{->}"4";"10"};
{\ar@{->}"7";"5"};
\end{xy}
\]
where $\sigma_0:D_0\ra \P_2$ is the blow up map and $s_0$ is the canonical section of $\O_{D_0}(L)$. Using the bottom row of this diagram one concludes that the Hilbert polynomial of $\sigma_0^*(\cal F)$ with respect to the sheaf $\O_{D_0}(1,1)$ is $6m+1$.
The restriction of~\refeq{eq:resolution back} to $D_0$ becomes the exact sequence
\[
2\O_{L}(-1)\xra{\varrho_{D_0}} \sigma_0^*(\cal F)\xra{\theta_{D_0}} \cal E|_{D_0}\ra 0.
\]
Therefore, we conclude that the Hilbert polynomial of the kernel of $\varrho_{D_0}$ is zero, hence $\varrho_{D_0}$ is injective.

By the shape of the original matrix $A$ we conclude that $\tilde A_{D_0}=B\cdot\pmat{
u_1& u_1y_1+u_2y_2\\
u_2&u_1z_1+u_2z_2
}$, where $B$ is identified with its corresponding $2\times 2$ matrix.
If $B=0$, then $\tilde A_{D_0}=0$ and $\varrho_{D_0}=0$, which contradicts the injectivity of $\varrho_{D_0}$.

Assume that the rank of $B$ is $1$ and that
\[
 B=\pmat{\lambda \alpha&\lambda \beta\\\mu \alpha&\mu \beta},\quad  (\alpha, \beta)\neq (0,0), (\lambda, \mu)\neq (0,0).
\]
In this case the kernel of $B$ is isomorphic to $\O_{D}(-2H)$ and is generated by the matrix $\pmat{\mu&-\lambda}$.
Then $\tilde A_{D_0}=\smat{\lambda u&\lambda q\\\mu u&\mu q}$ for $u=\alpha u_1+\beta u_2$ and
$q=\alpha(u_1y_1+u_2y_2)+\beta(u_1z_1+u_2z_2)$. Then the kernel of $\tilde A_{D_0}$
 is generated by $\pmat{\mu&-\lambda}$ and is isomorphic to $\O_{D}(-H-F)$.
As $\varrho_{D_0}$ is injective, we conclude that the kernels of $B$ and $\tilde A_{D_0}$ must be isomorphic, which is impossible since $H\not\sim F$.
\end{proof}
Since $B$ is an isomorphism, using the action of $\GL_2(\k)$ on the upper row of diagram~\refeq{eq:lifting res}, we may assume that $B=\id$.
Then $\tilde A$ can be written as
\[
\tilde A=
\pmat{
u_1&u_1y_1+u_2y_2\\
u_2&u_1z_1+u_2z_2
}+
\pmat{\xi_0&\xi_{00}x_0\\\eta_0&\eta_{00}x_0}u_0
\]
 using decomposition~\refeq{eq:decomposition of x_i}.

\par
\noindent
3)
By the resolution of $\sigma^*\cal F$ we obtain the following exact diagram
\begin{equation}\label{eq:diagram inverse}
\begin{xy}
(0,25)*+{}="31";
(15,25)*+{}="32";
(40,25)*+{ 0}="33";
(75,25)*+{}="34";
(100,25)*+{0}="35";
(115,25)*+{}="36";
(0,15)*+{}="21";
(15,15)*+{}="22";
(40,17)*+{ 2\O_{D_1}(-L)}="23";
(75,15)*+{}="24";
(100,15)*+{2\O_{D_1}(-L)}="25";
(115,15)*+{}="26";
(-8,0)*+{0}="11";
(8,0)*+{2\O_{D_1}(-L)}="12";
(40,0)*+{ 2\O_{D}(-2H)}="13";
(75,0)*+{ \O_{D}(-H)\oplus\O_{D}}="14";
(100,0)*+{\sigma^*\cal F }="15";
(115,0)*+{ 0}="16";
(0,-15)*+{}="01";
(15,-15)*+{0}="02";
(40,-15)*+{\cal K}="03";
(75,-15)*+{ \O_{D}(-H)\oplus\O_{D}}="04";
(100,-15)*+{\cal E}="05";
(115,-15)*+{ 0}="06";
(0,-25)*+{}="-11";
(15,-25)*+{}="-12";
(40,-25)*+{2\O_{D_1}(-L)}="-13";
(75,-25)*+{}="-14";
(100,-25)*+{0}="-15";
(115,-25)*+{}="-16";
(40,-35)*+{0}="-23";
{\ar@{->}"11";"12"};
{\ar@{->}_-{\smat{u_0&0\\0&u_0}}"12";"13"};
{\ar@{->}^-{\sigma^*A}"13";"14"};
{\ar@{->}^-{\pi}"14";"15"};
{\ar@{->}"15";"16"};
{\ar@{->}"02";"03"};
{\ar@{->}^-{\Phi}"03";"04"};
{\ar@{->}^-{\pi'}"04";"05"};
{\ar@{->}"05";"06"};
{\ar@{->}"33";"23"};
{\ar@{->}^{\smat{u_0&0\\0&u_0}}"23";"13"};
{\ar@{->}"13";"03"};
{\ar@{->}"03";"-13"};
{\ar@{->}"-13";"-23"};
{\ar@{->}"35";"25"};
{\ar@{->}^{\varrho}"25";"15"};
{\ar@{->}^{\theta}"15";"05"};
{\ar@{->}"05";"-15"};
{\ar@{=}"14";"04"};
{\ar@{=}"12";"23"};
\end{xy}
\end{equation}
with $\cal K$ the kernel of $\theta\circ \pi$.

By the construction of $\tilde A$, $\theta \circ \pi\circ \tilde A=0$, and
 hence $\tilde A$ factorizes uniquely through $\cal K$, with the commutative diagram
\[
\begin{xy}
(-10,15)*+{2\O_{D}(-H-F)}="1";
(20,0)*+{\O_{D}(-H)\oplus \O_{D}. }="2";
(-10,0)*+{\cal K}="3";
{\ar@{>->}^-{\Phi}"3";"2"};
{\ar@{>->}^{\tilde A}"1";"2"};
{\ar@{-->}_{\exists!\ \alpha}"1";"3"};
\end{xy}
\]
Because $\Phi$ and $\tilde A$ are injective, $\alpha$ is injective as well. On the other hand, $\cal K$ and $2\O_{D}(-H-F)$ have the same Hilbert polynomial with respect to $\O_{D}(H+F)$, namely $2(2m^2-m)$. Hence $\alpha$ is an isomorphism. Replacing
$\cal K$ by  $2\O_{D}(-H-F)$ yields the diagram of the proposition. Note that $\tilde A= \Phi(A, B)$ for
$B=\smat{\xi_0x_0&\xi_{00}x_0^2\\\eta_0x_0&\eta_{00}x_0^2}$.
This completes the proof of
Proposition~\ref{pr:all R-bundles are degenerations}.
\end{proof}

\subsection*{Proof of Proposition~\ref{pr:R properties}.}
1) Follows directly from~Proposition~\ref{pr:all R-bundles are degenerations}.

2) Restricting the resolution of $\cal E$ to $D_1$ one obtains the required
statement.

3) Restricting to $D_0$ we see that $\cal E|_{D_0}$ is the cokernel of
$\pmat{
u_1& u_1y_1+u_2y_2\\
u_2&u_1z_1+u_2z_2
}.$

Because $\smat{u_1\\u_2}$ is surjective, one obtains that $\cal E|_{D_0}$ is
the structure sheaf of the curve given by the determinant of this matrix.

4) Since
 $\Phi$ in~\refeq{eq:diagram} is of the form $\smat{l_1&q_1\\l_2&q_2}$ such
that
  $l_1, l_2\in \Gamma(D, \O_D(F))$ are linearly independent  with a single
common zero $q\in D_1\minus L$, there is the exact sequence
\[
0\ra \O_{D}(-2F-H)\xra{\smat{l_2&-l_1}}
2\O_{D}(-F-H)\xra{\smat{l_1\\l_2}}
\O_{D}(-H)\ra \k_q\ra 0.
\]
Splitting this sequence into short exact sequences one obtains the exact diagram
\[
\begin{xy}
(-45,25)*+{0}="-12";
(0,25)*+{0}="02";
(40,25)*+{}="12";
(-65,15)*+{0}="-21";
(-45,15)*+{\O_{D}(-2F-H)}="-11";
(0,15)*+{\O_{D}}="01";
(30,15)*+{\O_{C}}="11";
(40,15)*+{0}="21";
(-65,0)*+{0}="-20";
(-45,0)*+{2\O_{D}(-H-F)}="-10";
(0,0)*+{\O_{D}(-H)\oplus \O_{D}}="00";
(30,0)*+{\cal E}="10";
(40,0)*+{0}="20";
(-65,-15)*+{0}="-2-1";
(-45,-15)*+{\cal A}="-1-1";
(0,-15)*+{\O_{D}(-H)}="0-1";
(30,-15)*+{\k_q}="1-1";
(40,-15)*+{0}="2-1";
(-45,-25)*+{0}="-1-2";
(0,-25)*+{0}="0-2";
(30,-25)*+{}="1-2";
{\ar@{->}"-21";"-11"};
{\ar@{->}^-{-(l_1q_2-l_2q_1)}"-11";"01"};
{\ar@{->}"01";"11"};
{\ar@{->}"11";"21"};
{\ar@{->}"-20";"-10"};
{\ar@{->}^{\smat{l_1&q_1\\l_2& q_2}}"-10";"00"};
{\ar@{->}"00";"10"};
{\ar@{->}"10";"20"};
{\ar@{->}"-2-1";"-1-1"};
{\ar@{->}"-1-1";"0-1"};
{\ar@{->}"0-1";"1-1"};
{\ar@{->}"1-1";"2-1"};
{\ar@{->}"-12";"-11"};
{\ar@{->}_{\smat{l_2&-l_1}}"-11";"-10"};
{\ar@{->}"-10";"-1-1"};
{\ar@{->}"-1-1";"-1-2"};
{\ar@{->}"02";"01"};
{\ar@{->}^{\smat{0&1}}"01";"00"};
{\ar@{->}^{\smat{1\\0}}"00";"0-1"};
{\ar@{->}"0-1";"0-2"};
{\ar@{->}_{\smat{l_1\\l_2}}"-10";"0-1"};
\end{xy}
\]
and then, by the snake lemma, the required extension
\[
0\ra \O_C\ra \cal E\ra \k_q\ra 0.
\]
If this extension is trivial, then ${\sigma}_*\k_q\iso \k_p$ must be a direct summand of the $(3m+1)$-sheaf
${\sigma}_*\cal E\iso \cal F$ (cf. Remark~\ref{rem:push forward}).
This contradicts the stability of $\cal F$.

\begin{rem}
The point $q$ from Proposition~\ref{pr:R properties}, 4) arises from a  flat degeneration as above as follows. The points $p(t)$ of the fibres $\fcal F_t$ of the flat family $\fcal F$ form a section $S$ of $T\times \P_2$ over $T$ passing through $(0, \point{1, 0, 0})$. After blowing up its proper transform meats $D$ in the point $q$.
\end{rem}
\begin{rem}\label{rem:diagram for 3m+1}
Let $\cal F$ be a $(3m+1)$-sheaf as in the proof of Proposition~\ref{pr:all R-bundles are degenerations}. Then  already $\cal F$   induces an exact commutative diagram of type~\refeq{eq:diagram}. However  in this case
$\Phi=\smat{
u_1&u_1y_1+u_2y_2\\
u_2&u_1z_1+u_2z_2
}$ and its  cokernel is a one-dimensional sheaf $\tilde{\cal F}$ such that $\tilde{\cal F}|_{D_1}$ is a singular $(2m+2)$-sheaf on $D_1\iso \P_2$ whereas  $\tilde{\cal F}|_{D_0}$ is the structure sheaf of the supporting curve $C_0$ as before.
\end{rem}

\begin{rem}\label{rem:resolutions and maps}
Note that because of the vanishing of the relevant $Ext$-groups
every morphism between sheaves on $D(p)$ with resolutions of
the type~\refeq{eq:resolution new}
can be uniquely lifted to a morphism of the corresponding resolutions.
In particular this holds for $R$-bundles.
\end{rem}


\section{Classification result (main result)}
\label{section:main result}

In this section we are going to prove theorem \ref{tr:equivalence of sheaves}.
First of all note that the relation ``to be equivalent''
defined in Definition~\ref{df:equivalence for new sheaves} is in fact
an equivalence relation on the set of $R$-bundles associated to a
given singular $(3m+1)$-sheaf.

\begin{proof}[Proof of Theorem~\ref{tr:equivalence of sheaves}]
Since the morphism $X\ra M$ induces an isomorphism $N_{[\cal F]}\iso N_A$, where  $N_A=N_A(X_8)=T_A(X)/T_A(X_8)$ is the normal space to $X_8$ at the point $A$, it is enough to show that every two $R$-bundles $\cal E_1=\cal E(A, B_1)$ and $\cal E_2=\cal E(A, B_2)$ on $D(p)$ as in~\refeq{eq:E(A, B)}
are equivalent if and only if
$B_1$ and $B_2$ represent the same point in $\P N_A$.

Without loss of generality one can assume that $p=\point{1,0,0}$ and that $A$
is of the form~\refeq{eq:A special}. We can write $A$ as in~\refeq{eq:A special coeff}.
Adding a suitable multiple of the first column of $A$ to the second column
(this gives us an affine automorphism of $X$)
we can assume
without loss of generality that the coefficients $A_{01}$,
$A_{11}$, and $A_{12}$ are zero.

Let $\cal E_1=\cal E(A, B_1)$ and $\cal E_2=\cal E(A, B_2)$ be two equivalent
$R$-bundles, then the sheaves $\cal E_1$ and $\cal E_2$ possess locally free
resolutions of type~\refeq{eq:Phi(A, B)},
they are cokernels of $\Phi_1=\Phi(A, B_1)$ and $\Phi_2=\Phi(A, B_2)$ respectively.

Equivalence of $\cal E_1$  and $\cal E_2$ means that there exists an isomorphism $\phi:D(p)\ra D(p)$ identical on
$D_0(p)$ such that there is an isomorphism $\cal E_2\xra{\xi} \phi^*(\cal E_1)$.
By Remark~\ref{rem:resolutions and maps} $\xi$ can be uniquely lifted to a morphism of resolutions
\begin{equation}\label{eq:morph}
\begin{xy}
(0,0)*+{\O_{D(p)}(-H)\oplus \O_{D(p)}}="1";
(30,0)*+{\cal E_2}="2";
(0,-15)*+{\O_{D(p)}(-H)\oplus \O_{D(p)}}="3";
(30,-15)*+{\phi^*(\cal E_1)}="4";
{\ar@{->}"1";"2"};
{\ar@{->}"3";"4"};
{\ar@{->}^{
\smat{\bar a& \bar b \\0&\bar d}
}"1";"3"};
{\ar@{->}^{\xi}"2";"4"};
(-45,0)*+{2\O_{D(p)}(-H-F)}="5";
(-45,-15)*+{2\O_{D(p)}(-H-F)}="6";
{\ar@{->}^{\smat{a&b\\c&d}}"5";"6"};
{\ar@{->}^-{\phi^*(\Phi_1)}"6";"3"};
{\ar@{->}^-{\Phi_2}"5";"1"};
(-70, 0)*+{0}="7";
(-70,-15)*+{0}="8";
{\ar@{->}"8";"6"};
{\ar@{->}"7";"5"};
(45,0)*+{0}="9";
(45,-15)*+{0.}="10";
{\ar@{->}"2";"9"};
{\ar@{->}"4";"10"};
\end{xy}
\end{equation}
Note that from the uniqueness of the lifting it follows that isomorphisms
between $R$-bundles lift to automorphisms of  $\O_{D(p)}(-H)\oplus
\O_{D(p)}$ and thus the induced endomorphisms of $2\O_{D(p)}(-H-F)$ are also
automorphisms in this case. Therefore,  both matrices $\smat{a&b\\c&d}$ and $\smat{\bar a& \bar b \\0&\bar d}$ are invertible.

Straightforward verifications (cf.~\cite{MyGermanDiss}) show that for some $\mu\in \k^*$ the matrix
 $B_2-\mu B_1$ satisfies the tangent equations~\refeq{eq:tangent equations},
 i.~e., $B_2-\mu B_1\in T_A(X_8)$. So $B_1$
and $B_2$ represent the same element in $\P N_A$.

Let now $B_1$ and $B_2$
 be two equivalent normal vectors at $A\in X_8$.
Let
$\Phi_1=\Phi(A, B_1)$ and $\Phi_2=\Phi(A, B_2)$
be the matrices defining as in~\refeq{eq:Phi(A, B)} the sheaves
 $\cal E_1$ and $\cal E_2$ respectively.
Since $B_1$ and $B_2$ define the same point in $\P N_A$,
it follows that
\[
B_2-\alpha \cdot B_1\in T_A(X_8)
\]
for some $\alpha\in \k^*$.

Let
\[
B_1=\pmat{
\xi_0x_0+\xi_1x_1+\xi_2x_2&\xi_{00}x_0^2+\dots+\xi_{22}x_2^2\\
\eta_0x_0+\eta_1x_1+\eta_2x_2&\eta_{00}x_0^2+\dots+\eta_{22}x_2^2
}
\]
and
\[
 B_2=
\pmat{
\mu_0x_0+\mu_1x_1+\mu_2x_2&\mu_{00}x_0^2+\dots+\mu_{22}x_2^2\\
\nu_0x_0+\nu_1x_1+\nu_2x_2&\nu_{00}x_0^2+\dots+\nu_{22}x_2^2
}.
\]
Take \[
\beta=\mu_0-\xi_0\alpha,\quad \gamma=\nu_0-\eta_0 \alpha,
\]
and let
\begin{equation}\label{eq:L automorphism}
\phi_1=\smat{\alpha&\beta&\gamma\\
0&1&0\\
0&0&1
}:\P_2\ra \P_2,\quad \point{u_0, u_1, u_2}\mto \langle(u_0, u_1, u_2)
\smat{\alpha&\beta&\gamma\\
0&1&0\\
0&0&1
}\rangle.
\end{equation}
Note that the automorphisms of the form~\refeq{eq:L automorphism} are exactly
the automorphisms of $D_1\iso\P_2$ acting identically on $L$.
Consider now such an automorphism
$\phi:D(p)\ra D(p)$
with $\phi|_{D_1}=\phi_1$ and $\phi|_{D_0}=\id_{D_0}$.
Using the tangent equations~\refeq{eq:tangent equations} one checks that  $\phi^*(\Phi_1)=\Phi_2$. Therefore, there is an
isomorphism $\phi^*(\cal E_1)\iso \cal E_2$,
which means that the sheaves $\cal E_1$ and $\cal E_2$ are equivalent.
This completes the proof.
\end{proof}
\begin{rem} When one considers isomorphism classes of the $R$-bundles
 one finds that they depend on the point $q$ in Proposition
\ref{pr:R properties}. However, the isomorphism classes of the restrictions
of the $R$-bundles to the plane $D_1(p)$ don't depend on that point, as well
as the equivalence classes of the $R$-bundles.
\end{rem}

\section{Examples}
\label{section:examples}
Let us fix some $A\in X_8$ and let us consider the $R$-bundles $\cal E(A, B)$.
The curve $C_0$ is uniquely defined by the matrix $A$. The curve $C_1$ depends on $B$.
Let us fix $C_1$. Then by Proposition~\ref{pr:R properties}, 4), points of $C_1\minus L$ parameterize the isomorphism classes of $R$-bundles with the support $C_0\cup C_1$. To be more precise: there is a one-to-one correspondence between the isomorphism classes of $R$-bundles supported on $C_0\cup C_1$ and non-singular points of $C_1\minus L$.

It may however happen that there are two different points of the curve $C_1$ that define the same equivalence class of $R$-bundles.
This is the case when the curve $C_1$ has a non-trivial stabilizer under the action of the automorphisms of $D_1$ identical on $L$ (automorphisms from~\refeq{eq:L automorphism}). There is a natural action of this stabilizer on $C_1\minus L$ (and also on its non-singular locus) and the orbits of this action are clearly in one to one correspondence with the equivalence classes of $R$-bundles $\cal E(A, B)$ with support $C_0\cup C_1$.

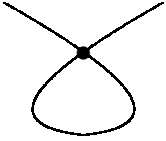
\begin{wrapfigure}[6]{l}{5cm}
\begin{picture}(110,70)(-55,-15)
\qbezier(-30,50)(50,5)(0,0)
\qbezier(30,50)(-50,5)(0,0)
\put(0,31){\circle*{5}}
\put(-50,-15){\small$x_1^2x_0-x_2(x_2+x_0)=0$}
\end{picture}
\end{wrapfigure}
\subsection*{Generic case: $C_0\cap L$ consists of two points.}
Let us fix the matrix $A=\pmat{x_1&x_2(x_0+x_2)\\x_2&x_1x_0}$. Let $C$ be a curve in $\P_2$ given by the determinant of this matrix. This is an irreducible cubic curve with an ordinary double point singularity.

 Then for
all directions $B$ the curve $C_0$ is given by the determinant of the matrix
\(
\pmat{u_1&u_2(x_0+x_2)\\u_2&u_1x_0}.
\)
The intersection of $C_0$ with the
line $L$ is given by the equation $u_1^2-u_2^2=0$ and consists of two
points, say  $q_1$ and $q_2$.
For a direction
\[
B=\pmat{
\xi_0x_0+\xi_1x_1+\xi_2x_2&\xi_{00}x_0^2+\dots+\xi_{22}x_2^2\\
\eta_0x_0+\eta_1x_1+\eta_2x_2&\eta_{00}x_0^2+\dots+\eta_{22}x_2^2
}
\]
 the restriction of $\cal E(A, B)$ to $D_1$ is given as the cokernel of the
matrix
 \[
 \pmat{u_1+\xi_0u_0&u_2+\xi_{00}u_0\\
u_2+\eta_0u_0&u_1+\eta_{00}u_0}.
\]
Its  support   is then the conic in $D_1$ through the points $q_1$ and $q_2$ given by the
determinant of this matrix.
\begin{rem}[comparison with~\cite{Seshadri}]
One sees that  $C_0$ is a normalization of $C$ ($C_0$ is a proper transform of $C$).

If the conic $C_1$ is smooth, then it is isomorphic to $\P_1$ and thus
the support of  an $R$-bundle in this situation is a curve of type $X_1$ (see~\cite{Seshadri}, pp. 212--213).

If $C_1$ is singular, then it is just a union of two lines and thus the whole support $C_0\cup C_1$ is a curve of type $X_2$.
\end{rem}
Let us fix $C_1$.

\subsubsection*{Smooth curve $C_1$.}
As already noticed, the points of $C_1\minus L$ are in one-to-one correspondence with the isomorphism classes of the $R$-bundles supported on $C_0\cup C_1$.
One sees in this case that the stabilizer group of $C_1$ consists of two elements.
The non-trivial one is the central symmetry with respect to the intersection point of the tangent lines to $C_1$ at $q_1$ and $q_2$ respectively (axes of $C_1$).

The orbits of the action of the stabilizer group on $C_1\minus L$ consist clearly of two points, the orbit space may be identified with the curve $C_1\minus L$ itself, the quotient map being induced by the $2:1$ self-covering $\k^*\ra \k^*$, $t\mto t^2$, via an isomorphism $C_1\minus L\iso \k^*$.
\subsubsection*{Singular curve $C_1$.}
In this case $C_1$ is a union of two lines. The stabilizer group of $C_1$ is a group isomorphic to $k^*$ and there are only two orbits in the non-singular locus  of $C_1\minus L$ each being one of the components of $C_1$ without their intersection point and without the points on $L$. Each orbit is isomorphic to $\k^*$.
\begin{figure}
\begin{center}
\begin{picture}
(140,150)(-25,-65)
\put(0,0){\line(4,1){80}}
\put(0,0){\line(0,1){60}}
\put(0,0){\line(0,-1){60}}
\put(80,20){\line(0,-1){60}}
\put(80,20){\line(0,1){60}}
\put(0,60){\line(4,1){80}}
\put(0,-60){\line(4,1){80}}
\put(0,0){\line(1,-1){40}}
\put(80,20){\line(1,-1){40}}
\put(40, -40){\line(4, 1){80}}
\put(0, 0){\line(-1, 1){25}}
\put(80, 20){\line(-1, 1){25}}
\put(-25, 25){\line(4, 1){80}}
\qbezier(5,55)(50,-80)(75,75)
\put(26,6.5){\circle*{5}}
\put(62,15.5){\circle*{5}}
\qbezier(40,-30)(-11,70)(105,-10)
\put(0,0){\line(-4,-1){10}}
\put(80,20){\line(4,1){10}}
\put(82, 25){\small $L$}
\put(10,50){\small ${C_0}$}
\put(100,-20){\small$C_1$}
\put(100,0){\small$D_1$}
\put(40,76){\small$D_0$}
\end{picture}
\quad
\begin{picture}
(140,150)(-25,-65)
\put(0,0){\line(4,1){80}}
\put(0,0){\line(0,1){60}}
\put(0,0){\line(0,-1){60}}
\put(80,20){\line(0,-1){60}}
\put(80,20){\line(0,1){60}}
\put(0,60){\line(4,1){80}}
\put(0,-60){\line(4,1){80}}
\put(0,0){\line(1,-1){40}}
\put(80,20){\line(1,-1){40}}
\put(40, -40){\line(4, 1){80}}
\put(0, 0){\line(-1, 1){25}}
\put(80, 20){\line(-1, 1){25}}
\put(-25, 25){\line(4, 1){80}}

\qbezier(5,55)(50,-80)(75,75)
\put(26,6.5){\circle*{5}}
\put(62,15.5){\circle*{5}}
\put(26,6.5){\line(5,-2){50}}
\put(26,6.5){\line(-5,2){20}}
\put(62,15.5){\line(1,-5){7}}
\put(62,15.5){\line(-1,5){4}}

\put(0,0){\line(-4,-1){10}}
\put(80,20){\line(4,1){10}}
\put(82, 25){\small $L$}
\put(10,50){\small ${C_0}$}
\put(55,-18){\small$C_1$}
\put(100,0){\small$D_1$}
\put(40,76){\small$D_0$}
\end{picture}
\end{center}
\caption{Support of $R$-bundles in the generic case.}\label{fig:Support generic}
\end{figure}
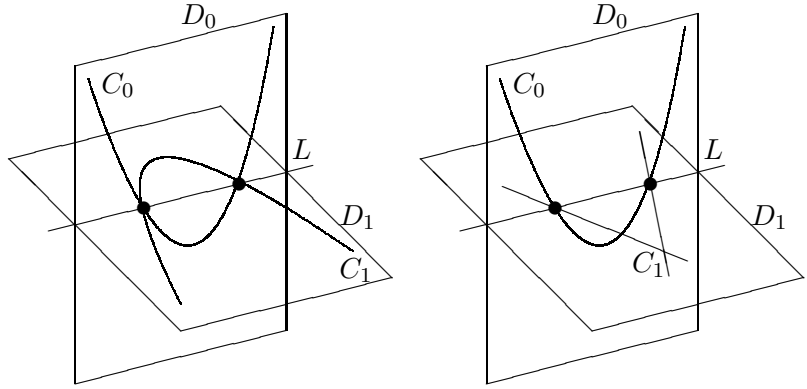

\medskip

Another cases with two points in the intersection $C_0\cap L$ are the following.

\subsubsection*{Three lines with simple intersections.}
\begin{figure}
\begin{center}
\begin{picture}
(140,150)(-25,-65)
\put(0,0){\line(4,1){80}}
\put(0,0){\line(0,1){60}}
\put(0,0){\line(0,-1){60}}
\put(80,20){\line(0,-1){60}}
\put(80,20){\line(0,1){60}}
\put(0,60){\line(4,1){80}}
\put(0,-60){\line(4,1){80}}
\put(0,0){\line(1,-1){40}}
\put(80,20){\line(1,-1){40}}
\put(40, -40){\line(4, 1){80}}
\put(0, 0){\line(-1, 1){25}}
\put(80, 20){\line(-1, 1){25}}
\put(-25, 25){\line(4, 1){80}}

\put(26,6.5){\line(-3,5){25}}
\put(26,6.5){\line(3,-5){15}}

\put(62,15.5){\line(1,5){12}}
\put(62,15.5){\line(-1,-5){5}}
\put(2,35){\line(5,2){75}}

\put(26,6.5){\circle*{5}}
\put(62,15.5){\circle*{5}}
\qbezier(40,-30)(-11,70)(105,-10)
\put(0,0){\line(-4,-1){10}}
\put(80,20){\line(4,1){10}}
\put(82, 25){\small $L$}
\put(10,50){\small ${C_0}$}
\put(100,-20){\small$C_1$}
\put(100,0){\small$D_1$}
\put(40,76){\small$D_0$}
\end{picture}
\quad
\begin{picture}
(140,150)(-25,-65)
\put(0,0){\line(4,1){80}}
\put(0,0){\line(0,1){60}}
\put(0,0){\line(0,-1){60}}
\put(80,20){\line(0,-1){60}}
\put(80,20){\line(0,1){60}}
\put(0,60){\line(4,1){80}}
\put(0,-60){\line(4,1){80}}
\put(0,0){\line(1,-1){40}}
\put(80,20){\line(1,-1){40}}
\put(40, -40){\line(4, 1){80}}
\put(0, 0){\line(-1, 1){25}}
\put(80, 20){\line(-1, 1){25}}
\put(-25, 25){\line(4, 1){80}}
\put(0,0){\line(-4,-1){10}}
\put(80,20){\line(4,1){10}}

\put(26,6.5){\line(-3,5){25}}
\put(26,6.5){\line(3,-5){15}}
\put(62,15.5){\line(1,5){12}}
\put(62,15.5){\line(-1,-5){5}}
\put(2,35){\line(5,2){75}}

\put(26,6.5){\circle*{5}}
\put(62,15.5){\circle*{5}}

\put(26,6.5){\line(5,-2){50}}
\put(26,6.5){\line(-5,2){20}}
\put(62,15.5){\line(1,-5){7}}
\put(62,15.5){\line(-1,5){4}}

\put(82, 25){\small $L$}
\put(10,50){\small ${C_0}$}
\put(55,-18){\small$C_1$}
\put(100,0){\small$D_1$}
\put(40,76){\small$D_0$}
\end{picture}
\end{center}
\caption{Support of $R$-bundles for three lines with simple intersections.}
\label{fig:3 lines}
\end{figure}
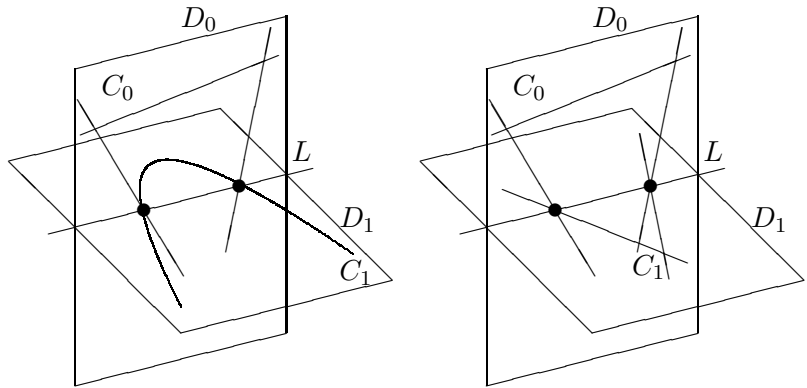

As a typical example one can choose the matrix $A=\pmat{x_1&0\\x_2&x_2x_0}$.
We present the pictures for different types of $C_1$ in Figure~\ref{fig:3 lines}.
\subsubsection*{Transversal intersection of a line with a smooth conic.}
A typical example comes from the matrix $A=\pmat{x_1&x_0x_1\\x_2&x_1^2}$. The corresponding pictures can be seen on Figure~\ref{fig: line and conic}
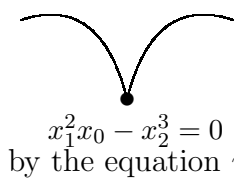
\begin{figure}
\begin{center}
\begin{picture}
(140,150)(-25,-65)
\put(0,0){\line(4,1){80}}
\put(0,0){\line(-4,-1){10}}
\put(80,20){\line(4,1){10}}
\put(0,0){\line(0,1){60}}
\put(0,0){\line(0,-1){60}}
\put(80,20){\line(0,-1){60}}
\put(80,20){\line(0,1){60}}
\put(0,60){\line(4,1){80}}
\put(0,-60){\line(4,1){80}}
\put(0,0){\line(1,-1){40}}
\put(80,20){\line(1,-1){40}}
\put(40, -40){\line(4, 1){80}}
\put(0, 0){\line(-1, 1){25}}
\put(80, 20){\line(-1, 1){25}}
\put(-25, 25){\line(4, 1){80}}

\put(62,15.5){\line(1,5){12}}
\put(62,15.5){\line(-1,-5){5}}

\qbezier(40,-10)(-20,55)(77,55)
\put(70,55){\circle*{2}}

\qbezier(40,-30)(-11,70)(105,-10)

\put(26,6.5){\circle*{5}}
\put(62,15.5){\circle*{5}}

\put(82, 25){\small $L$}
\put(10,45){\small ${C_0}$}
\put(100,-20){\small$C_1$}
\put(100,0){\small$D_1$}
\put(40,76){\small$D_0$}
\end{picture}
\quad
\begin{picture}
(140,150)(-25,-65)
\put(0,0){\line(4,1){80}}
\put(0,0){\line(-4,-1){10}}
\put(80,20){\line(4,1){10}}
\put(0,0){\line(0,1){60}}
\put(0,0){\line(0,-1){60}}
\put(80,20){\line(0,-1){60}}
\put(80,20){\line(0,1){60}}
\put(0,60){\line(4,1){80}}
\put(0,-60){\line(4,1){80}}
\put(0,0){\line(1,-1){40}}
\put(80,20){\line(1,-1){40}}
\put(40, -40){\line(4, 1){80}}
\put(0, 0){\line(-1, 1){25}}
\put(80, 20){\line(-1, 1){25}}
\put(-25, 25){\line(4, 1){80}}

\put(62,15.5){\line(1,5){12}}
\put(62,15.5){\line(-1,-5){5}}

\qbezier(40,-10)(-20,55)(77,55)
\put(70,55){\circle*{2}}

\put(26,6.5){\circle*{5}}
\put(62,15.5){\circle*{5}}

\put(26,6.5){\line(5,-2){50}}
\put(26,6.5){\line(-5,2){20}}
\put(62,15.5){\line(1,-5){7}}
\put(62,15.5){\line(-1,5){4}}

\put(82, 25){\small $L$}
\put(10,45){\small ${C_0}$}
\put(100,-20){\small$C_1$}
\put(100,0){\small$D_1$}
\put(40,76){\small$D_0$}
\end{picture}
\end{center}
\caption{Support of $R$-bundles for a simple intersection of a line and a smooth conic.}\label{fig: line and conic}
\end{figure}
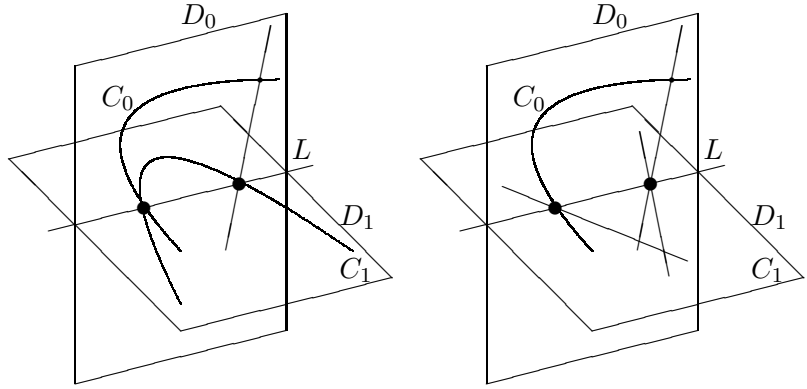

\begin{wrapfigure}[3]{l}{5cm}
\begin{picture}(80,50)(-45,-23)
\qbezier(0,0)(-10,40)(-40,30)
\qbezier(0,0)(10,40)(40,30)
\put(0,0){\circle*{5}}
\put(-30,-15){\small$x_1^2x_0-x_2^3=0$}
\end{picture}
\end{wrapfigure}

\subsection*{$C_0\cap L$ consists of a single point.}
Let us fix the matrix $A=\pmat{x_1&x_2^2\\x_2&x_1x_0}$. The curve $C$ is in this case a cuspidal cubic curve.
 The intersection of $C_0$ with the line
$L$ is given by the equation $u_1^2=0$ and consists of a single point.
Let us fix some $C_1$.

\subsubsection*{Smooth curve $C_1$.}
One sees that the stabilizer of $C_1$ is trivial in this case. Hence the points of  $C_1\minus L$ are in one-to-one correspondence with the isomorphism and simultaneously with the equivalence classes of $R$-bundles supported on $C_0\cup C_1$.

\subsubsection*{Singular curve $C_1$.}
In this case the stabilizer of $C_1$ is a $1$-dimensional group acting transitively on $C_1\minus L$.
\begin{figure}
\begin{center}
\begin{picture}(140,150)(-25,-65)
\put(0,0){\line(4,1){80}}
\put(0,0){\line(0,1){60}}
\put(0,0){\line(0,-1){60}}
\put(80,20){\line(0,-1){60}}
\put(80,20){\line(0,1){60}}
\put(0,60){\line(4,1){80}}
\put(0,-60){\line(4,1){80}}
\put(0,0){\line(1,-1){40}}
\put(80,20){\line(1,-1){40}}
\put(40, -40){\line(4, 1){80}}
\put(0, 0){\line(-1, 1){25}}
\put(80, 20){\line(-1, 1){25}}
\put(-25, 25){\line(4, 1){80}}
\put(0,0){\line(-4,-1){10}}
\put(80,20){\line(4,1){10}}

\qbezier(5,55)(40,-45)(75,75)

\qbezier(40,-30)(5,40)(105,-10)

\put(40,10){\circle*{5}}

\put(100,-20){\small$C_1$}
\put(10,50){\small ${C_0}$}
\put(82, 25){\small $L$}
\put(100,0){\small$D_1$}
\put(40,76){\small$D_0$}
\end{picture}
\quad
\begin{picture}(140,150)(-25,-65)
\put(0,0){\line(4,1){80}}
\put(0,0){\line(0,1){60}}
\put(0,0){\line(0,-1){60}}
\put(80,20){\line(0,-1){60}}
\put(80,20){\line(0,1){60}}
\put(0,60){\line(4,1){80}}
\put(0,-60){\line(4,1){80}}
\put(0,0){\line(1,-1){40}}
\put(80,20){\line(1,-1){40}}
\put(40, -40){\line(4, 1){80}}
\put(0, 0){\line(-1, 1){25}}
\put(80, 20){\line(-1, 1){25}}
\put(-25, 25){\line(4, 1){80}}

\qbezier(5,55)(40,-45)(75,75)
\put(40,10){\circle*{5}}

\put(40,10){\line(1,-5){8}}
\put(40,10){\line(-1,5){3}}

\put(40,10){\line(5,-4){35}}
\put(40,10){\line(-5,4){10}}

\put(0,0){\line(-4,-1){10}}
\put(80,20){\line(4,1){10}}
\put(50,-25){\small$C_1$}
\put(10,50){\small ${C_0}$}
\put(82, 25){\small $L$}
\put(100,0){\small$D_1$}
\put(40,76){\small$D_0$}
\end{picture}
\end{center}
\caption{}\label{fig:one point cusp}
\end{figure}

\medskip

There are also the following similar cases.
\subsubsection*{Tangent intersection of a line with a smooth conic.}
An example for this case is provided by the matrix $A=\pmat{x_1&x_0x_2\\x_2&x_1x_2}$. The pictures are given in Figure~\ref{fig:tangent}.
\begin{figure}
\begin{center}
\begin{picture}(140,150)(-25,-65)
\put(0,0){\line(4,1){80}}
\put(0,0){\line(0,1){60}}
\put(0,0){\line(0,-1){60}}
\put(80,20){\line(0,-1){60}}
\put(80,20){\line(0,1){60}}
\put(0,60){\line(4,1){80}}
\put(0,-60){\line(4,1){80}}
\put(0,0){\line(1,-1){40}}
\put(80,20){\line(1,-1){40}}
\put(40, -40){\line(4, 1){80}}
\put(0, 0){\line(-1, 1){25}}
\put(80, 20){\line(-1, 1){25}}
\put(-25, 25){\line(4, 1){80}}
\put(0,0){\line(-4,-1){10}}
\put(80,20){\line(4,1){10}}

\put(40,10){\line(-1,5){10}}
\put(40,10){\line(1,-5){4}}

\qbezier(75,-35)(2,-15)(75,65)

\qbezier(40,-30)(5,40)(105,-10)

\put(40,10){\circle*{5}}

\put(100,-20){\small$C_1$}
\put(36,32){\small ${C_0}$}
\put(82, 25){\small $L$}
\put(100,0){\small$D_1$}
\put(40,76){\small$D_0$}
\end{picture}
\quad
\begin{picture}(140,150)(-25,-65)
\put(0,0){\line(4,1){80}}
\put(0,0){\line(0,1){60}}
\put(0,0){\line(0,-1){60}}
\put(80,20){\line(0,-1){60}}
\put(80,20){\line(0,1){60}}
\put(0,60){\line(4,1){80}}
\put(0,-60){\line(4,1){80}}
\put(0,0){\line(1,-1){40}}
\put(80,20){\line(1,-1){40}}
\put(40, -40){\line(4, 1){80}}
\put(0, 0){\line(-1, 1){25}}
\put(80, 20){\line(-1, 1){25}}
\put(-25, 25){\line(4, 1){80}}
\put(0,0){\line(-4,-1){10}}
\put(80,20){\line(4,1){10}}

\put(40,10){\line(-1,5){10}}
\put(40,10){\line(1,-5){4}}

\qbezier(75,-35)(2,-15)(75,65)


\put(40,10){\line(4,-1){60}}
\put(40,10){\line(-4,1){10}}

\put(40,10){\line(1,-1){35}}
\put(40,10){\line(-1,1){10}}

\put(40,10){\circle*{5}}

\put(95,-15){\small$C_1$}
\put(36,32){\small ${C_0}$}
\put(82, 25){\small $L$}
\put(100,0){\small$D_1$}
\put(40,76){\small$D_0$}
\end{picture}
\end{center}
\caption{}\label{fig:tangent}
\end{figure}
\subsubsection*{Point on a double line.}
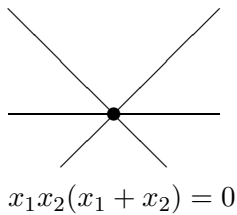
\begin{figure}
\begin{center}
\begin{picture}(140,150)(-25,-65)
\put(0,0){\line(4,1){80}}
\put(0,0){\line(0,1){60}}
\put(0,0){\line(0,-1){60}}
\put(80,20){\line(0,-1){60}}
\put(80,20){\line(0,1){60}}
\put(0,60){\line(4,1){80}}
\put(0,-60){\line(4,1){80}}
\put(0,0){\line(1,-1){40}}
\put(80,20){\line(1,-1){40}}
\put(40, -40){\line(4, 1){80}}
\put(0, 0){\line(-1, 1){25}}
\put(80, 20){\line(-1, 1){25}}
\put(-25, 25){\line(4, 1){80}}
\put(0,0){\line(-4,-1){10}}
\put(80,20){\line(4,1){10}}

\put(39.3,10){\line(-1,5){10}}
\put(39.3,10){\line(1,-5){4}}
\put(40.7,10){\line(-1,5){10}}
\put(40.7,10){\line(1,-5){4}}

\put(5,45){\line(5,1){70}}

\qbezier(40,-30)(5,40)(105,-10)

\put(40,10){\circle*{5}}

\put(100,-20){\small$C_1$}
\put(17,53){\small ${C_0}$}
\put(82, 25){\small $L$}
\put(100,0){\small$D_1$}
\put(40,76){\small$D_0$}
\end{picture}
\begin{picture}
(140,150)(-25,-65)
\put(0,0){\line(4,1){80}}
\put(0,0){\line(0,1){60}}
\put(0,0){\line(0,-1){60}}
\put(80,20){\line(0,-1){60}}
\put(80,20){\line(0,1){60}}
\put(0,60){\line(4,1){80}}
\put(0,-60){\line(4,1){80}}
\put(0,0){\line(1,-1){40}}
\put(80,20){\line(1,-1){40}}
\put(40, -40){\line(4, 1){80}}
\put(0, 0){\line(-1, 1){25}}
\put(80, 20){\line(-1, 1){25}}
\put(-25, 25){\line(4, 1){80}}
\put(0,0){\line(-4,-1){10}}
\put(80,20){\line(4,1){10}}

\put(39.3,10){\line(-1,5){10}}
\put(39.3,10){\line(1,-5){4}}
\put(40.7,10){\line(-1,5){10}}
\put(40.7,10){\line(1,-5){4}}

\put(5,45){\line(5,1){70}}


\put(40,10){\line(4,-1){60}}
\put(40,10){\line(-4,1){10}}

\put(40,10){\line(1,-1){35}}
\put(40,10){\line(-1,1){10}}

\put(40,10){\circle*{5}}

\put(95,-15){\small$C_1$}
\put(17,53){\small ${C_0}$}
\put(82, 25){\small $L$}
\put(100,0){\small$D_1$}
\put(40,76){\small$D_0$}
\end{picture}
\end{center}
\caption{}\label{fig:only double line}
\end{figure}

We can consider the matrix $A=\pmat{x_1&0\\x_2&x_0x_1}$. The pictures are given in Figure~\ref{fig:only double line}.

\begin{wrapfigure}[5]{l}{4cm}
\begin{picture}(85,75)(-40,-35)
\put(0,0){\line(1,1){40}}
\put(0,0){\line(-1,-1){20}}
\put(0,0){\line(-1,1){40}}
\put(0,0){\line(1,-1){20}}
\put(0,0){\line(1,0){40}}
\put(0,0){\line(-1,0){40}}
\put(0,0){\circle*{5}}
\put(-40,-35){\small$x_1x_2(x_1+x_2)=0$}
\end{picture}
\end{wrapfigure}
\subsection*{$C_0\cap L$ is the whole line $L$.}
We start here from the matrix $A=\pmat{x_1&0\\x_2&x_2(x_1+x_2)}$, $A=\pmat{x_1&0\\x_2&x_2^2}$ or $A=\pmat{x_1&0\\x_2&x_1^2}$.
Then the curve $C$ consists of three different lines that all intersect in the
same point, of a line and a double line, or of a triple line respectively.

In this case $C_1$ is the union the line $L$ with another
line $L_1$, which intersects $L$ at a point $p_B$ that depends on the
direction $B$.
It is clear that the stabilizer of $C_1$ acts transitively on $L_1$, hence there is only one equivalence class of $R$-bundles with fixed $C_1$ in this case. Two directions  $B$ and
$B'$ with different intersection points $p_B$ and $p_{B'}$ define
non-equivalent sheaves because all the allowed automorphism  are
identities on $L$. See Figure~\ref{fig: with double line} for the corresponding pictures.

\begin{figure}
\begin{center}
\begin{picture}
(140,150)(-25,-65)
\put(0,0){\line(4,1){80}}
\put(0,0){\line(-4,-1){10}}
\put(80,20){\line(4,1){10}}
\put(0,0){\line(0,1){60}}
\put(0,0){\line(0,-1){60}}
\put(80,20){\line(0,-1){60}}
\put(80,20){\line(0,1){60}}
\put(0,60){\line(4,1){80}}
\put(0,-60){\line(4,1){80}}
\put(0,0){\line(1,-1){40}}
\put(80,20){\line(1,-1){40}}
\put(40, -40){\line(4, 1){80}}
\put(0, 0){\line(-1, 1){25}}
\put(80, 20){\line(-1, 1){25}}
\put(-25, 25){\line(4, 1){80}}

\put(20,5){\line(-1,3){15}}
\put(20,5){\line(1,-3){10}}
\put(40,10){\line(-1,5){10}}
\put(40,10){\line(1,-5){7}}
\put(64,16){\line(1,6){8}}
\put(64,16){\line(-1,-6){7}}

\put(20,5){\circle*{5}}

\put(40,10){\circle*{5}}

\put(64,16){\circle*{5}}

\put(52,13){\circle*{5}}
\put(52,13){\line(-1,1){20}}
\put(52,13){\line(1,-1){35}}

\put(82, 25){\small $L$}
\put(100,0){\small$D_1$}
\put(51,18){\tiny$p_B$}
\put(40,76){\small$D_0$}
\put(90,-20){\small$L_1$}
\end{picture}
\quad
\begin{picture}
(140,150)(-25,-65)
\put(0,0){\line(4,1){80}}
\put(0,0){\line(-4,-1){10}}
\put(80,20){\line(4,1){10}}
\put(0,0){\line(0,1){60}}
\put(0,0){\line(0,-1){60}}
\put(80,20){\line(0,-1){60}}
\put(80,20){\line(0,1){60}}
\put(0,60){\line(4,1){80}}
\put(0,-60){\line(4,1){80}}
\put(0,0){\line(1,-1){40}}
\put(80,20){\line(1,-1){40}}
\put(40, -40){\line(4, 1){80}}
\put(0, 0){\line(-1, 1){25}}
\put(80, 20){\line(-1, 1){25}}
\put(-25, 25){\line(4, 1){80}}

\put(20,5){\line(-1,3){15}}
\put(20,5){\line(1,-3){10}}
\put(63.3,16){\line(1,6){8}}
\put(63.3,16){\line(-1,-6){7}}
\put(64.7,16){\line(1,6){8}}
\put(64.7,16){\line(-1,-6){7}}

\put(20,5){\circle*{5}}


\put(64,16){\circle*{5}}

\put(52,13){\circle*{5}}
\put(52,13){\line(-1,1){20}}
\put(52,13){\line(1,-1){35}}

\put(82, 25){\small $L$}
\put(100,0){\small$D_1$}
\put(51,18){\tiny$p_B$}
\put(40,76){\small$D_0$}
\put(90,-20){\small$L_1$}
\end{picture}
\quad
\begin{picture}
(140,150)(-25,-65)
\put(0,0){\line(4,1){80}}
\put(0,0){\line(-4,-1){10}}
\put(80,20){\line(4,1){10}}
\put(0,0){\line(0,1){60}}
\put(0,0){\line(0,-1){60}}
\put(80,20){\line(0,-1){60}}
\put(80,20){\line(0,1){60}}
\put(0,60){\line(4,1){80}}
\put(0,-60){\line(4,1){80}}
\put(0,0){\line(1,-1){40}}
\put(80,20){\line(1,-1){40}}
\put(40, -40){\line(4, 1){80}}
\put(0, 0){\line(-1, 1){25}}
\put(80, 20){\line(-1, 1){25}}
\put(-25, 25){\line(4, 1){80}}

\put(18.8,5){\line(-1,3){15}}
\put(18.8,5){\line(1,-3){10}}
\put(20,5){\line(-1,3){15}}
\put(20,5){\line(1,-3){10}}
\put(21.2,5){\line(-1,3){15}}
\put(21.2,5){\line(1,-3){10}}


\put(20,5){\circle*{5}}



\put(52,13){\circle*{5}}
\put(52,13){\line(-1,1){20}}
\put(52,13){\line(1,-1){35}}

\put(82, 25){\small $L$}
\put(100,0){\small$D_1$}
\put(51,18){\tiny$p_B$}
\put(40,76){\small$D_0$}
\put(90,-20){\small$L_1$}

\end{picture}
\end{center}
\caption{}\label{fig: with double line}
\end{figure}

\section{Universal families of $R$-bundles}
We are going to construct flat families parameterizing all $R$-bundles as
well as all non-singular $(3m+1)$-sheaves.

\label{section:universal family}
\subsection*{Family over $\tilde X$, construction of the space}
Let $\tilde X\xra{\alpha} X$ be the blow up of $X$ along $X_8$ and
let $H=X\times h$ for a line $h\subset \P_2$. Then there is
 the universal $(3m+1)$-sheaf $\fcal U$  on $X\times \P_2$  given by the locally free resolution (cf.~\cite{Freiermuth}, 6.1)
\begin{equation}\label{eq:universsal 3m+1}
0\ra 2\O_{X\times\P_2}(-2H)\xra{A_X} \O_{X\times\P_2}(-H)\oplus
\O_{X\times\P_2}\ra \fcal U\ra 0,
\end{equation}
where  $A_X$ is the universal matrix: $A_X|_{\{A\}\times\P_2}= A$ for all
$A\in X$.

By pulling back we
obtain the family $\bar{ \fcal U}\defeq (\alpha\times \id)^*\fcal U$ of
$(3m+1)$-sheaves over $\tilde X$.  Let $S_8=\Sing\fcal U$ be the closed
subvariety of $X\times \P_2$ where $\fcal U$
is not locally free, i.~e., the zero locus of $A_X$.
Since $X$ is the parameter space of the universal cubic,
the subvariety $S_8$ is a section of $X\times\P_2$ over $X_8$ and so isomorphic
to $X_8.$ In particular $S_8$ is smooth of codimension 3.
Then $\tilde S_8\defeq(\alpha\times \id)^{-1}(S_8)$ is the set of
points in $\tilde X\times \P_2$ where the sheaf $\bar{\fcal U}$ is not
locally free.
$\tilde S_8$ is isomorphic to the exceptional divisor
$\tilde X_8=\alpha^{-1}(X_8)$ of the blow up $\tilde X\xra{\alpha} X$,
in particular $\tilde S_8$ is smooth.

Let
\(
\tau:Y\ra \tilde X\times \P_2
\)
be the blowing up of $\tilde X\times \P_2$ along $\tilde S_8$ and
let $\tilde D_1$ denote the exceptional divisor of
 $\tau$. Let $\tilde D_0$ be the proper transform of
$\tilde X_8\times \P_2$. It is isomorphic to the blow up of
$\tilde X_8\times \P_2$ along $\tilde S_8$. Then
$\tilde D=\tilde D_0\cup \tilde D_1$ is the preimage of
 $\tilde X_8\times \P_2$.
\[
\begin{xy}
(0,0)*+{\tilde D_1}="1";
(20,0)*+{Y}="2";
(0,-15)*+{\tilde S_8}="3";
(20,-15)*+{\tilde X\times \P_2}="4";
{\ar@{->}^-{}"1";"3"};
{\ar@{>->}^-{}"1";"2"};
{\ar@{>->}^-{}"3";"4"};
{\ar@{->}^-{\tau}"2";"4"};
\end{xy}
\hskip 2cm
\begin{xy}
(0,0)*+{\tilde X\times \P_2}="1";
(25,0)*+{X\times\P_2}="2";
(0,-15)*+{\tilde X}="3";
(25,-15)*+{X}="4";
(-15,-5)*+{\tilde S_8}="5";
(10,-5)*+{S_8}="6";
(-15,-20)*+{\tilde X_8}="7";
(10,-20)*+{X_8}="8";
{\ar@{->}^{\alpha\times \id}"1";"2"};
{\ar@{->}^-(0.7){\alpha}|(0.4){\hole}"3";"4"};
{\ar@{->}^-{pr_1}|(0.35){\hole}"1";"3"};
{\ar@{->}^-{pr_1}"2";"4"};
{\ar@{->}"5";"6"};
{\ar@{->}"7";"8"};
{\ar@{->}^-{\iso}"5";"7"};
{\ar@{->}_-{\iso}"6";"8"};
{\ar@{>->}^-{}"5";"1"};
{\ar@{>->}^-{}"6";"2"};
{\ar@{>->}^-{}"7";"3"};
{\ar@{>->}^-{}"8";"4"};
\end{xy}
\]
A fibre $Y_x$ of the morphism $Y\to\tilde X\times\P_2\to\tilde X$ for
a point $x\in\tilde X_8$ is a surface $D(p)$ with $p\in\tilde S_8$
corresponding to $x$. Moreover, the composed morphism
$Y\to\tilde X$ is flat. Hence the situation is analogous to that of the blow
up $Z$ in section \ref{section:1-fam}. There is also an embedding of $Y$
analogous to that of $Z$:

\begin{rem}\label{rem:remp2bundle}
The ideal sheaf $\ki$ of $\tilde S_8$ is the quotient of a decomposable
rank-$3$ vector bundle as follows. Let
$x_1, x_2\in\Gamma\ko_{\tilde X\times\P_2}(H)$ be the independent linear
entries of the lifted matrix $A$. Then $\ki$ is generated by
$x_1, x_2$ and the lifted equation of $\tilde X_8$, giving rise to a
surjection
\[
\E:=\ko_{\tilde X}(-\tilde X_8)\boxtimes\ko_{\P_2}\oplus 2\ko_{\tilde X\times\P_2}(-H)\to\ki.
\]
It follows that $Y=\P(\ki)$ is embedded in the $\P_2$-bundle $\P(\E)$
with local equations similar to those of $Z.$
\end{rem}

\subsection*{Family over $\tilde X$, construction of the sheaf}
Pulling back sequence~\refeq{eq:universsal 3m+1} along the morphism
\[
Y\xra{\tau}\tilde X\times \P_2\xra{\alpha\times \id} X\times \P_2,
\]
we obtain the sequence
\[
0\ra 2\O_{Y}(-2H)\xra{} \O_{Y}(-H)\oplus \O_{Y}\ra
\tau^*(\bar {\fcal U})\ra 0,
\]
where $H$ also denotes the pull-back of $\tilde{X}\times h.$
This remains exact because the sheaf $O_{Y}(-2H)$ is locally free and,
therefore, has no torsion.
Similar to diagram~\refeq{eq:main construction diagram} one obtains the commutative diagram with exact rows and columns
\begin{equation}\label{eq:family over tilde X}
\begin{xy}
(105,20)*+{0}="35";
(35,10)*+{0}="23";
(105,10)*+{\cal C}="25";
(10,0)*+{0}="12";
(35,0)*+{ 2\O_{Y}(-2H)}="13";
(75,0)*+{ \O_{Y}(-H)\oplus \O_{Y}}="14";
(105,0)*+{\tau^*\bar{\fcal U}}="15";
(115,0)*+{0,}="16";
(10,-15)*+{0}="02";
(35,-15)*+{2\O_{Y}(-2H+\tilde D_1)}="03";
(75,-15)*+{ \O_{Y}(-H)\oplus \O_{Y}}="04";
(105,-15)*+{\tilde{\fcal U}}="05";
(115,-15)*+{0}="06";
(35,-25)*+{\cal C}="-13";
(105,-25)*+{0}="-15";
(35,-35)*+{0}="-23";
{\ar@{->}"12";"13"};
{\ar@{->}^-{}"13";"14"};
{\ar@{->}"14";"15"};
{\ar@{->}"15";"16"};
{\ar@{->}"02";"03"};
{\ar@{->}^-{}"03";"04"};
{\ar@{->}"04";"05"};
{\ar@{->}"05";"06"};
{\ar@{->}"23";"13"};
{\ar@{->}^{}"13";"03"};
{\ar@{->}"03";"-13"};
{\ar@{->}"-13";"-23"};
{\ar@{->}"35";"25"};
{\ar@{->}"25";"15"};
{\ar@{->}"15";"05"};
{\ar@{->}"05";"-15"};
{\ar@{=}"14";"04"};
\end{xy}
\end{equation}
by multiplying with the canonical section of $\O_Y(\tilde{D}_1)$, with
$\cal C=\O_{\tilde D_1}\ten\O_Y(-2H+\tilde D_1)$, where $\tilde{\fcal U}$
is the quotient.

\begin{pr}\label{pr:flatn}
1) The sheaf\; $\tilde{\fcal U}$ is flat over $\tilde X$
and the fibres of\;  $\tilde{\fcal U}$ are either non-singular
$(3m+1)$-sheaves or $R$-bundles on some $D(p).$

2) Any non-singular $(3m+1)$-sheaf and any $R$-bundle on some $D(p)$
is equivalent to a fibre of\; $\tilde{\fcal U}$.
\end{pr}
\begin{proof}Consider $l_B$ as in~\refeq{eq:l_B}, an embedding of an open set
of $\k$ in $X$ along a
normal direction $B\in \k^{18}$ such that $0\in T$ is the only point in $T$
with the image in $X_8$. Then from the  universal property of blow-ups it
follows that
$l_B$ uniquely factorizes through $\tilde X\xra{\alpha} X$, i. e.,
there exists the commutative diagram
\[
\xymatrix{
&\tilde X\ar[d]^{\alpha}\\
U\ar[r] ^{l_B}\ar@{.>}^{\exists!\ \tilde l_B}[ur]&X.
}
\]
Then using again the universal property of blow-ups we obtain the commutative
diagram with cartesian squares
\begin{equation}\label{eq:diagram 1 and tilde X}
\begin{xy}
(0,15)*+{Z}="0";
(30,15)*+{Y}="00";
(0,0)*+{T\times \P_2}="1";
(30,0)*+{\tilde X\times \P_2}="2";
(60,0)*+{X\times \P_2}="3";
(0,-15)*+{T}="4";
(30,-15)*+{\tilde X}="5";
(60,-15)*+{X.}="6";
{\ar@{->}^{\tilde l_B\times \id}"1";"2"};
{\ar@{->}^{\alpha\times \id}"2";"3"};
{\ar@{->}^{\tilde l_B}"4";"5"};
{\ar@{->}^{\alpha}"5";"6"};
{\ar@{->}"1";"4"};
{\ar@{->}"2";"5"};
{\ar@{->}"3";"6"};
{\ar@{.>}^{\exists!\ L_B}"0";"00"};
{\ar@{->}^{\sigma}"0";"1"};
{\ar@{->}^{\tau}"00";"2"};
\end{xy}
\end{equation}
Comparing diagram~\refeq{eq:diagram 1 and tilde X} with the construction
of the blow up $Z$ in Section~\ref{section:R-bundles}, one finds that the
restriction of the exceptional divisor and the proper transform of
$\tilde X_8\times\P_2$ to the image of $T$ in $\tilde X$ is the one of $Z$.
It follows that the restriction of diagram~\refeq{eq:family over tilde X}
to this image is isomorphic to diagram~\refeq{eq:main construction diagram}.
It follows that the sheaf $\tilde{\fcal U}$ is flat over $\tilde X$ and
locally free on its support.

In particular, the construction of $R$-bundles
(diagram~\refeq{eq:main construction diagram}) is obtained by pulling back
diagram~\refeq{eq:family over tilde X} from $Y$ along $L_B$ as
in~\refeq{eq:diagram 1 and tilde X}. In particular this means that the fibres
of the sheaf $\tilde{\fcal U}$ are either non-singular $(3m+1)$-sheaves
or $R$-bundles on some $D(p)$.
By the above construction of $l_B$ every non-singular $(3m+1)$-sheaf
(up to isomorphy) as well as every $R$-bundle (up to equivalence)
appear as fibres of $\tilde{\fcal U}$.
\end{proof}

\begin{rem}
Because $\tilde{\fcal U}$ is flat over  $\tilde X$ it is easy to show that
the sheaf $\cal C$ is the ``relative torsion'' of $\tau^*\bar{\fcal U}$
over $\tilde{S}_8$. The same holds for diagram
\refeq{eq:main construction diagram}.
\end{rem}

\subsection*{The universal family over $\tilde M$}
Let $Y_M=\Bl_{\tilde \nu(\tilde S_8)}(\tilde M\times \P_2)$ and let $Y_M\xra{\tau_M} \tilde M\times \P_2$ be the corresponding morphism. Let $\tilde D_M$ be the exceptional divisor.
Then by the universal property of blow ups there exists a unique
morphism  $Y\xra{\xi} Y_M$ such that the diagram
\begin{equation}\label{eq:Y_M}
\begin{split}
\begin{xy}
(0,0)*+{Y}="1";
(30,0)*+{Y_M}="2";
(0,-15)*+{\tilde X\times\P_2}="3";
(30,-15)*+{\tilde M\times\P_2}="4";
{\ar@{->}^{\tilde \nu\times \id} "3";"4"};
{\ar@{->}^{\tau} "1";"3"};
{\ar@{->}^{\tau_M} "2";"4"};
{\ar@{->}^{\xi} "1";"2"};
\end{xy}
\end{split}
\end{equation}
commutes.

\begin{pr}\label{pr:family over tilde M}
There exists a flat sheaf $\tilde{\fcal V}$ on $Y_M$ such that
$\tilde{\fcal U}$ is a pull back of $\tilde{\fcal V}$ along $\xi$ up to a
twist by a pull back of a line bundle on $\tilde X$.
\end{pr}

\begin{proof} Let $\fcal V$ be the universal $(3m+1)$-family on $M\times \P_2$.
Let $p_1:M\times \P_2\ra M$ and $p_2:M\times \P_2\ra \P_2$ be the canonical projections. Then there exists a relative Beilinson resolution (cf.~\cite{Okonek-Schneider-Spindler})
\begin{equation}\label{eq:relative Beilinson M}
0\ra p_1^*\cal A_2\ten p_2^*\O_{\P_2}(-2)\xra{}
(p_1^*\cal A_1\ten p_2^*\O_{\P_2}(-1))
\oplus
(p_1^*\cal A_0\ten p_2^*\O_{\P_2})\ra
 \fcal V\ra 0,
\end{equation}
where $\cal A_2=R^1{p_1}_{*}(\fcal V\ten p_2^*\Omega_{\P_2}^2(2))$,
$\cal A_1=R^1{p_1}_{*}(\fcal V\ten p_2^*\Omega_{\P_2}^1(1))$, and $\cal A_0=R^0{p_1}_{*}(\fcal V)$  are locally free sheaves on $S$ of rank $2$, $1$, and $1$ respectively.

Consider the blow up $\alpha_M:\tilde M\ra M$.
Let $\bar{\fcal V}$ be the pull-back of $\fcal V$ to $\tilde M\times \P_2$ along $\alpha_M\times \id$. Pulling back once more along $\tau_M$ one obtains the exact sequence
\[
0\ra
\tau_M^*(\cal A_2\sqten\O_{\P_2}(-2))
\xra{}
\tau_M^*(\cal A_1\sqten \O_{\P_2}(-1)
\oplus
\cal A_0\sqten \O_{\P_2})
\ra
\tau_M^*(\bar {\fcal V})\ra 0.
\]
Using the canonical section of $\O_{Y_M}(\tilde D_M)$, similarly to the construction of~\refeq{eq:family over tilde X} we obtain the commutative diagram with exact rows
\begin{equation}\label{eq:family over tilde M}
\begin{xy}
(0,7)*+{0}="12";
(35,7)*+{\tau_M^*(\cal A_2\sqten\O_{\P_2}(-2))}="13";
(98,7)*+{\tau_M^*(\cal A_1\sqten \O_{\P_2}(-1)
\oplus
\cal A_0\sqten \O_{\P_2})
}="14";
(135,7)*+{\tau_M^*\bar{\fcal V}}="15";
(145,7)*+{0,}="16";
(0,-8)*+{0}="02";
(35,-8)*+{\tau_M^*(\cal A_2\sqten\O_{\P_2}(-2))\ten \O_{Y_M}(\tilde D_M)
)}="03";
(98,-8)*+{\tau_M^*(\cal A_1\sqten \O_{\P_2}(-1)
\oplus
\cal A_0\sqten \O_{\P_2})
}="04";
(135,-8)*+{\tilde{\fcal V}}="05";
(145,-8)*+{0}="06";
{\ar@{->}"12";"13"};
{\ar@{->}^-{}"13";"14"};
{\ar@{->}"14";"15"};
{\ar@{->}"15";"16"};
{\ar@{->}"02";"03"};
{\ar@{->}^-{}"03";"04"};
{\ar@{->}"04";"05"};
{\ar@{->}"05";"06"};
{\ar@{->}^{}"13";"03"};
{\ar@{->}"15";"05"};
{\ar@{=}"14";"04"};
\end{xy}
\end{equation}
Note that since the morphism  $X\xra{\nu}M$ is flat, the pull-back
of~\refeq{eq:relative Beilinson M} along $\nu\times \id$ is the relative
Beilinson resolution of $\fcal U$ up to a twist by $p_1^*\cal L$ for some
line bundle $\cal L$ on $X$, i.~e., resolution~\refeq{eq:universsal 3m+1}
tensorized by $p_1^*\cal L$. Therefore, the commutativity of~\refeq{eq:Y_M}
implies that the pull-back of~\refeq{eq:family over tilde M} along $\xi$
coincides with~\refeq{eq:family over tilde X} twisted by the pull back of
the line bundle $\nu^*\cal L$ on $\tilde X$.
This completes the proof of Proposition~\ref{pr:family over tilde M}.
\end{proof}

\begin{cor} The family $\tilde{\fcal V}$ is universal in the following sense.
Each equivalence class of $R$-bundles has a unique representative
as the fibre of $\tilde{\fcal V}$ over its corresponding point in $\tilde M_8$
and each isomorphism class of non-singular $(3m+1)$-sheaves has a unique representative
as the fibre of $\tilde{\fcal V}$ over its corresponding point in
$\tilde{M}\minus\tilde M_8$.
Moreover, the sheaves on $Z$ as in Section~\ref{section:1-fam}
(one-dimensional deformations of $R$-bundles) are obtained as pull-backs
of $\tilde{\fcal V}$ via unique morphisms $T\to\tilde M$.
\end{cor}



\def\cprime{$'$} \def\cprime{$'$}

\end{document}